\documentclass{amsart}
\usepackage{times,mathrsfs,array,amssymb,epsfig,tikz}
\usetikzlibrary{petri}

\newcounter{lemma}
\newtheorem{Theorem}[lemma]{Theorem}
\newtheorem{Lemma}[lemma]{Lemma}
\newtheorem{Corollary}[lemma]{Corollary}
\newtheorem{Proposition}[lemma]{Proposition}

\theoremstyle{definition}
\newtheorem{Example}[lemma]{Example}

\newtheorem{Notation}[lemma]{Notation}

\newtheorem{Remark}[lemma]{Remark}

\def\C{\mathbb C}
\def\HH{\mathbb H}
\def\H{\mathfrak H}

\def\Q{\mathbb Q}
\def\R{\mathbb R}
\def\O{\mathcal O}

\def\Z{\mathbb Z}
\def\new{\mathrm{new}}
\def\mod{\  \mathrm{mod}\ }
\def\Im{\mathrm{Im\,}}

\def\JS#1#2{\left(\frac{#1}{#2}\right)}

\def\SL{\mathrm{SL}}
\def\M#1#2#3#4{\begin{pmatrix}#1&#2\\#3&#4\end{pmatrix}}
\def\SM#1#2#3#4{\left(\begin{smallmatrix}#1&#2\\#3&#4\end{smallmatrix}
  \right)}
\def\gen#1{\langle #1\rangle}
\def\FL#1{\left\lfloor #1\right\rfloor}

\begin{document}

\title[Automorphic forms on Shimura curves]{Schwarzian differential
  equations and Hecke eigenforms on Shimura curves}
\author{Yifan Yang}
\address{Department of Applied Mathematics, National Chiao Tung
  University and National Center for Theoretical Sciences, Hsinchu,
  Taiwan 300}  
\email{yfyang@math.nctu.edu.tw}
\date{\today}
\thanks{The author was partially supported by Grant
  99-2115-M-009-011-MY3 of the National Science Council, Taiwan (R.O.C.).}

\begin{abstract} Let $X$ be a Shimura curve of genus zero. In this
  paper, we first characterize the spaces of automorphic forms on
  $X$ in terms of Schwarzian differential equations. We then devise a
  method to compute Hecke operators on these spaces. An interesting
  by-product of our analysis is the evaluation
  $$
    _2F_1\left(\frac1{24},\frac7{24};\frac56;
    -\frac{2^{10}\cdot3^3\cdot5}{11^4}\right)=\sqrt6
    \sqrt[6]{\frac{11}{5^5}}
  $$
  and other similar identities.
\end{abstract}
\maketitle

\begin{section}{Introduction}

Let $K$ be a totally real number field and $R$ be its ring of
integers. Let $B$ be a quaternion algebra over $K$ that splits exactly
at one infinite place, i.e.,
$$
  B\otimes_\Q\R\simeq M(2,\R)\times\HH^{[K:\Q]-1},
$$
where $M(2,\R)$ is the algebra of $2\times 2$ matrices over $\R$ and
$\HH$ is Hamilton's quaternion algebra. Thus, up to conjugation, there
is a unique embedding $\iota:B\hookrightarrow M(2,\R)$
from $B$ into $M(2,\R)$. Given an order $\O$ of $B$, let $\O_1^\ast$
denote the group of elements of reduced norm $1$ inside $\O$. Then the
image $\Gamma(\O):=\iota(\O_1^\ast)$ of $\O_1^\ast$ under $\iota$ is a
subgroup of $\SL(2,\R)$ and the quotient space
$\Gamma(\O)\backslash\H$ is called the \emph{Shimura curve} associated
to $\O$, where $\H$ denotes the upper half-plane. For example, when
$K=\Q$, $B=M(2,\Q)$, and $\O=M(2,\Z)$, we have $\O_1^\ast=\SL(2,\Z)$
and the Shimura curve associated to $M(2,\Z)$ is simply the classical
modular curve $Y_0(1)$. Thus, Shimura curves are generalizations of
classical modular curves. From now on, the term \emph{Shimura curve}
will be reserved strictly for the case $B\neq M(2,\Q)$. The compact
Riemann surface associated to a discrete subgroup of $\SL(2,\R)$
commensurable to some $\O_1^\ast$ will also be called
a Shimura curve. (See \cite[Chapter 2]{Elkies} and \cite[Chapter
2]{Alsina-Bayer} for more detailed introduction to Shimura curves.)

When $B$ is an indefinite quaternion algebra over $\Q$, Shimura curves
have moduli-space interpretation similar to their classical
counterpart. Namely, they are moduli spaces of
principally polarized abelian surfaces with quaternionic
multiplication. (See \cite{Shimura-CM1}.) Also, many
theories and properties about classical modular curves can be extended
to the case of Shimura curves. For example, if, for a nonnegative even
integer $k$, we let $S_k(\Gamma(\O))$ or simply $S_k(\O)$ denote the
space of automorphic forms of weight $k$, i.e., the space consisting
of holomorphic functions $f:\H\to\C$ satisfying
$$
  f\left(\frac{a\tau+b}{c\tau+d}\right)=(c\tau+d)^kf(\tau)
$$
for all $\SM abcd\in\Gamma(\O)$ and all $\tau\in\H$, then, similar to
the case of classical modular curves, we can define a family of
mutually commuting self-adjoint linear operators on $S_k(\O)$,
called \emph{Hecke operators}. (Here self-adjointness is with
respect to the Petersson inner product defined by
$$
  \gen{f,g}=\int_{\Gamma(\O)\backslash\H}f(\tau)\overline{g(\tau)}
  y^k\frac{dxdy}{y^2}, \qquad \tau=x+iy,
$$
where the integration is taken over any fundamental domain of
$\Gamma(\O)$.) Then the space $S_k(\O)$ contains a basis consisting of
simultaneous eigenforms, called \emph{Hecke eigenforms}, for all Hecke
operators. In fact, according to the Jacquet-Langlands correspondence
\cite[Chapter 16]{Jacquet-Langlands}, each Hecke eigenform corresponds
to an irreducible automorphic representation of $\mathrm{GL}(2,\Q)$
that is an inner twist of a certain irreducible cuspidal
representation. In other words, for each Hecke eigenform on
$\Gamma(\O)$, there corresponds a Hecke eigenform on a certain
modular group with the same eigenvalues.

On the other hand, even though it is true that many theoretical
aspects of classical modular curves can be extended to the case of
Shimura curves, it is not true for explicit methods. The main obstacle
lies at the lack of cusps on Shimura curves. Namely, in the case of
classical modular curves, many problems about modular curves can be
answered using $q$-expansions (i.e., expansions with respect to a
local parameter at a \emph{cusp}) of modular forms or modular
functions involved, and there are many explicit methods for
constructing modular functions and modular forms and computing their
$q$-expansions. In fact, because the Fourier coefficients of a
normalized Hecke eigenform on congruence subgroups are identical with
the eigenvalues of Hecke operators, one can compute the $q$-expansions
of Hecke eigenforms without actually constructing them. However,
because there are no cusps on Shimura curves, any method for classical
modular curves that uses $q$-expansions cannot possibly be extended to
the case of Shimura curves. Moreover, as far as we know, eigenvalues
for Hecke operators on automorphic forms on Shimura curves do not say
anything about Taylor coefficients of automorphic forms.
Thus, it is both interesting and challenging to find explicit methods
for Shimura curves.

In this paper, for a Shimura curve of genus zero, we will first
characterize the spaces of automorphic forms in terms of Schwarzian
differential equations. (See Remark \ref{remark: Schwarzian DE} for
the definition of Schwarzian differential equations.) In other words,
spaces of automorphic forms will be represented using solutions of
certain differential equations. This makes explicit computation on
automorphic forms possible. For example, in the second half of the
paper, we will devise a method to compute Hecke operators and hence
determine Hecke eigenforms. Two examples will be worked out. As
by-products of our analysis, we find the following intriguing
evaluations
\begin{equation*}
\begin{split}
  _2F_1\left(\frac1{24},\frac7{24};\frac56;
  -\frac{2^{10}\cdot3^3\cdot5}{11^4}\right)&=\sqrt6
  \sqrt[6]{\frac{11}{5^5}}, \\
  {}_2F_1\left(\frac5{24},\frac{11}{24};\frac76;
  -\frac{2^{10}\cdot3^3\cdot5}{11^4}\right)&=
   \frac{\sqrt[3]2(1+\sqrt2)\sqrt[6]{11^5}}{20\sqrt3}
   \frac{\Gamma(7/6)\Gamma(13/24)\Gamma(19/24)}
   {\Gamma(5/6)\Gamma(17/24)\Gamma(23/24)},
\end{split}
\end{equation*}
and other similar identities, where $_2F_1(a,b;c;x)$ denotes the
$_2F_1$-hypergeometric function. (See Corollary \ref{corollary:
  hypergeometric evaluation} and Remark \ref{remark: hypergeometric
  evaluation} below.)

Since Schwarzian differential equations play a crucial role in our
approach, it is an important problem to determine the differential
equation associated to each Shimura curve of genus $0$. Some work has
already been done in this direction. In particular, in \cite{Tu}, Tu
determined the Schwarzian differential equations associated to certain
Eichler orders in quaternion algebras over $\Q$.

As of now, our method only works for Shimura curves of genus zero. We
hope to extend our method to Shimura curves of higher genus in the
future. Note that Sijsling \cite{Sijsling} has already considered the
cases of Shimura curves of genus one with exactly one elliptic point
and obtained differential equations associated to these curves. It
will be an interesting problem to combine his equations with our
approach to study Shimura curves of genus one.
\end{section}

\begin{section}{Schwarzian differential equations and automorphic forms
    on Shimura curves}
Let $X(\O)$ be a Shimura curve with the associated norm-one group
$\Gamma(\O)$. It is known since the nineteenth century that if $F(\tau)$
is a meromorphic automorphic form of weight $k$ (with some multiplier
system) and $t(\tau)$ is a nonconstant automorphic function on
$\Gamma(\O)$, then $F,\tau F,\ldots,\tau^kF$, as functions of $t$,
span the solution space of a $(k+1)$-st order linear ordinary
differential equation
$$
  \theta^{k+1}F+r_k(t)\theta^kF+\cdots r_0(t)F=0, \quad \theta=t\frac d{dt},
$$
with algebraic functions as coefficients $r_j(t)$. (See \cite[Theorem
5.1]{Stiller} and also \cite[Theorem 1]{Yang}.) In
this paper, we refer to this kind of differential equations as
\emph{automorphic differential equations}. If
the compact Riemann surface $\Gamma\backslash\H$ has
genus zero, which we assume from now on, and $t(\tau)$ is a generator
of the function field on $\Gamma(\O)$, which we call a
\emph{Hauptmodul} of $\Gamma(\O)$, then the coefficients $r_j(t)$ are
actually rational functions and all the singularities of the
differential equations are regular. In fact, if $F(\tau)$ is a
holomorphic automorphic form, then the singularities of the
differential equation are precisely the points where the function
$\tau\to t(\tau)$ fails to be locally one-to-one, that is, where the
values of $t$ correspond to elliptic points. If the number of elliptic
points is $3$, i.e., if $\Gamma(\O)$ is a triangle group, then it is a
classical fact that a second order ordinary differential equation with
exactly three regular singular points is completely determined by the
local exponents. In this way, one can write down an automorphic
differential equation associated to a triangle group without actually
finding an automorphic form and a Hauptmodul first. (In fact, such a 
differential equation must be the symmetric power of an algebraic
transformation of a $_2F_1$-hypergeometric function.) Then one can
study properties of Shimura curves using this differential
equation. This method has been used by \cite{Elkies,Voight-CM} to study CM
(complex multiplication) points on Shimura curves.

When a Shimura curve $X(\O)$ of genus zero has more than three
elliptic points, the determination of automorphic differential
equations is more complicated. In \cite{Elkies}, Elkies determined a
differential equation associated to the normalizer of a maximal order
in a quaternion algebra of discriminant $10$ over $\Q$ and then used
this differential equation to numerically compute the coordinates of
CM points on a Hauptmodul. In \cite{Bayer-Travesa}, Bayer and Travesa
obtained automorphic differential equations for the Shimura curve
associated to the maximal order in the indefinite quaternion algebra
of discriminant $6$ and its various Atkin-Lehner quotients.
In both cases, the determination uses geometry of Shimura curves.
For instance, in \cite{Elkies}, Elkies used the covering between two
Shimura curves. Other than these two isolated results, as far as we
know, there is no systematic attempt in literature to determine such
differential equations. (Note that in \cite{Elkies}, such an
automorphic differential equation is called a Schwarz
equation. However, in this paper, we will reserve the 
term \emph{Schwarz equation} for a certain normalized automorphic
differential equation because of its connection to the Schwarzian
derivative. See Proposition \ref{proposition: Schwarzian} below.)

To fully realize the potential of the method of automorphic
differential equations, we shall first normalize such differential
equations. The idea is that if $t(\tau)$ is an automorphic function on
$\Gamma(\O)$, then $t'(\tau)$ is a meromorphic automorphic form of
weight $2$ on $\Gamma(\O)$. Thus, $t'(\tau)^{1/2}$, as a function of $t$,
satisfies a second-order ordinary differential equation. This
differential equation can be regarded as a normal form for all
automorphic differential equations associated to $\Gamma(\O)$ because
it depends only on the chosen automorphic function $t(\tau)$. In fact,
there is a simple formula to convert a general automorphic differential
equation to such a differential equation satisfied by $t'(\tau)$ and
$t(\tau)$.

\begin{Proposition}[{\cite[Pages 99--102,290]{Ford}}]
  \label{proposition: Schwarzian}
  Let $X(\O)$ be a Shimura curve. Let $F(\tau)$ be an automorphic form
  of weight $1$ (with some multiplier system) and $t(\tau)$ be a
  nonconstant automorphic function on $\Gamma(\O)$. If the
  second-order differential equation satisfied by $F(\tau)$ and
  $t(\tau)$ is
  $$
    \theta^2F+r_1(t)\theta F+r_0(t)F=0, \qquad \theta=t\frac d{dt},
  $$
  then the differential equation satisfied by $t'(\tau)^{1/2}$ and
  $t(\tau)$ is
  \begin{equation} \label{equation: Schwarzian DE}
    \frac{d^2}{dt^2}G+Q(t)G=0,
  \end{equation}
  where
  $$
    Q(t)=\frac{1+4r_0-2t(dr_1/dt)-r_1^2}{4t^2}
  $$
  and satisfies
  $$
    2Q(t)t'(\tau)^2+\{t,\tau\}=0, \qquad
    \{t,\tau\}=\frac{t'''(\tau)}{t'(\tau)}
   -\frac32\left(\frac{t''(\tau)}{t'(\tau)}\right)^2.
  $$
\end{Proposition}

The reader who is unfamiliar with the relation between automorphic
forms and differential equations may find the proof of this
proposition given in \cite{Yang} easier to comprehend.

\begin{Remark} \label{remark: Schwarzian DE}
The function $\{t,\tau\}$ is classically known as the
\emph{Schwarzian derivative}. (See \cite[Chapter 10]{Hille}.) In our
setting, it is a meromorphic automorphic form of weight $4$ on
$\Gamma(\O)$. In view of its connection to the Schwarzian derivative,
we call the differential equation in \eqref{equation: Schwarzian DE}
satisfied by $t'(\tau)^{1/2}$ and $t(\tau)$ the \emph{Schwarzian
  differential equation} associated to $t$. Note that the function
$-\{t,\tau\}/2t'(\tau)^2$, up to a factor of $-4$, is called the
\emph{automorphic derivative} in \cite{Bayer-Travesa}. Here we will
follow the same terminology.
\end{Remark}

\begin{Notation} For a thrice-differentiable function $f$ of $z$,
  we let $D(f,z)$ denote the automorphic derivative
  $$
    D(f,z)=-\frac{\{f,z\}}{2f'(z)^2}, \qquad
    \{f,z\}=\frac{f'''(z)}{f'(z)}
   -\frac32\left(\frac{f''(z)}{f'(z)}\right)^2.
  $$
\end{Notation}

Now the upshot is that if $X(\O)$ has genus zero and $t(\tau)$ is a
Hauptmodul, then the analytic behavior of $t'(\tau)$ is simple to
describe and it is easy to express all (holomorphic) automorphic forms
with trivial character in terms of $t'(\tau)$. Here the dimension
formula in the following theorem is taken from \cite{Shimura-book}.

\begin{Theorem} \label{theorem: basis}
  Assume that a Shimura curve $X$ has genus zero with elliptic
  points $\tau_1,\ldots,\tau_r$ of order $e_1,\ldots,e_r$,
  respectively. Let $t(\tau)$ be a Hauptmodul of $X$ and set
  $a_i=t(\tau_i)$, $i=1,\ldots,r$. For a positive even integer $k\ge 4$,
  let
  $$
    d_k=\dim S_k(\O)=1-k+\sum_{j=1}^r\left\lfloor
    \frac k2\left(1-\frac1{e_j}\right)\right\rfloor.
  $$
  Then a basis for the space of automorphic forms of weight $k$ on $X$ is
  $$
    t'(\tau)^{k/2}t(\tau)^j\prod_{i=1,a_i\neq\infty}^r
    \left(t(\tau)-a_j\right)^{-\lfloor k(1-1/e_1)/2\rfloor}, \quad
    j=0,\ldots,d_k-1.
  $$
\end{Theorem}

\begin{proof} If $a_j\neq\infty$, then we have
  $$
    t(\tau)-a_j=C(\tau-\tau_j)^{e_j}+O\left((\tau-\tau_j)^{e_j+1}\right),
    \quad C\in\C,
  $$
  near $\tau_j$ since $\tau_j$ is assumed to be elliptic point of
  order $e_j$. The constant $C$ cannot be zero because $t-a_j$ is also a
  Hauptmodul and cannot have a zero of order greater than $1$ (as a
  function on $X$) at $a_j$. Then
  \begin{equation} \label{equation: whatever}
    t'(\tau)=Ce_j(\tau-\tau_j)^{e_j-1}+O\left((\tau-\tau_j)^{e_j}\right)
  \end{equation}
  near $\tau_j$. Thus, the function
  \begin{equation} \label{equation: proposition basis}
    t'(\tau)^{k/2}\prod_{i=1,a_i\neq\infty}^r
    (t(\tau)-a_i)^{-\lfloor k(1-1/e_1)/2\rfloor}
  \end{equation}
  is a (possibly meromorphic) automorphic form of weight $k$ that is
  holomorphic throughout $\H$ except for possibly the point where
  $t(\tau)=\infty$.

  Consider first the case $a_j\neq\infty$ for all $j$. Let $\tau_0$ be
  the point where $t=\infty$. Since $t$ is a Hauptmodul, $t$ must have
  a simple pole at $\tau_0$, that is
  $$
    t(\tau)=\frac C{\tau-\tau_0}+O(1), \qquad C\in\C.
  $$
  Then $t'(\tau)=-C/(\tau-\tau_0)^2+O(1)$ and the order of the
  function in \eqref{equation: proposition basis} at $\tau_0$ is
  $$
    A:=-k+\sum_{j=1}^r\left\lfloor
    \frac k2\left(1-\frac1{e_j}\right)\right\rfloor.
  $$
  If $\dim S_k(\O)\neq 0$, then $A\ge 0$ and thus the function in
  \eqref{equation: proposition basis} is holomorphic throughout $\H$.
  In fact, we can multiply the function in \eqref{equation: proposition
    basis} by a polynomial of degree $\le A$ in $t$ and still get a
  function holomorphic throughout $\H$. Since the dimension of the
  space of polynomials of degree $\le A$ is the same as $d_k$, we
  conclude that the functions in the statement of the theorem form
  a basis for $S_k(\O)$.

  Now assume that $a_j=\infty$ for some $j$, say, $a_1=\infty$. Again,
  because $t$ is assumed to be a Hauptmodul, we must have
  $$
    t(\tau)=\frac C{(\tau-\tau_1)^{e_1}}
   +O\left((\tau-\tau_1)^{e_1-1}\right), \qquad C\in\C,
  $$
  near $\tau_1$. Then the order of the function \eqref{equation:
    proposition basis}, as a function of $\tau$, at $\tau_1$ is
  $$
    B:=-\frac k2(e_1+1)+e_1\sum_{i=2}^r\left\lfloor
    \frac k2\left(1-\frac1{e_j}\right)\right\rfloor.
  $$
  We can multiply the function in \eqref{equation: proposition basis}
  by a polynomial in $t$ of degree not exceeding
  $$
    \left\lfloor\frac B{e_1}\right\rfloor
   =\left\lfloor-\frac k2\left(1+\frac1{e_1}\right)\right\rfloor
    +\sum_{i=2}^r\left\lfloor
    \frac k2\left(1-\frac1{e_j}\right)\right\rfloor
   =-k+\sum_{j=1}^r\left\lfloor
    \frac k2\left(1-\frac1{e_j}\right)\right\rfloor
  $$
  and still get a function holomorphic throughout $\H$. Again, we
  conclude that the functions in the statement form a basis for
  $S_k(\O)$. This completes the proof.
\end{proof}

The combination of Proposition \ref{proposition: Schwarzian} and
Theorem \ref{theorem: basis} gives us a concrete space of functions
that can be used to study properties of automorphic forms, provided
that a Schwarzian differential equation has been determined. In the
second half of the paper, we will provide a method to compute Hecke
operators using these functions.

In view of the importance of Schwarzian differential equations, here
we shall review some analytic properties of Schwarzian differential
equation. These informations will be helpful in determining the
differential equations.

\begin{Proposition} \label{proposition: D}
  Automorphic derivatives have the following properties.
  \begin{enumerate}
  \item $D((az+b)/(cz+d),z)=0$ for all $\SM abcd\in\mathrm{GL}(2,\C)$.
  \item $D(g\circ f,z)=D(g,f)+D(f,z)/(dg/df)^2$.
 \end{enumerate}
\end{Proposition}

\begin{proof} Both properties are well-known and can be verified
  directly.
\end{proof}



\begin{Proposition} \label{proposition: Q}
  Assume that $X(\O)$ has genus zero with elliptic
  points $\tau_1,\ldots,\tau_r$ of order $e_1,\ldots,e_r$,
  respectively. Let $t(\tau)$ be a Hauptmodul of $X(\O)$ and set
  $a_i=t(\tau_i)$, $i=1,\ldots,r$. Then the automorphic derivative
  $Q(t)$ in Proposition \ref{proposition: Schwarzian} is equal to
  $$
    Q(t)=\frac14\sum_{j=1,a_j\neq\infty}^r\frac{1-1/e_j^2}{(t-a_j)^2}
        +\sum_{j=1,a_j\neq\infty}^r\frac{B_j}{t-a_j}
  $$
  for some constants $B_j$. Moreover, if $a_j\neq\infty$ for all $j$,
  then the constants $B_j$ satisfy
  $$
  \sum_{j=1}^r B_j=
  \sum_{j=1}^r\left(a_jB_j+\frac14(1-1/e_j^2)\right)=
  \sum_{j=1}^r\left(a_j^2B_j+\frac12a_j(1-1/e_j^2)\right)=0.
  $$
  Also, if $a_r=\infty$, then $B_j$ satisfy
  $$
  \sum_{j=1}^{r-1}B_j=0, \qquad
  \sum_{j=1}^{r-1}\left(a_jB_j+\frac14(1-1/e_j^2)\right)=\frac14(1-1/e_r^2).
  $$
\end{Proposition}

In principle, the properties stated above were known for a long time.
However, in literature, usually the proposition is stated under the
assumption that $a_i$ are all real. The reason for this assumption is
that in the standard books on ordinary differential equations in the
complex domain, one usually starts from a second-order linear
differential equation and build an automorphic function from it. At
some point, one would need to employ the Schwarz reflection
principle in order to extend the domain of the definition of the
automorphic function to the whole upper half-plane. This is where the
assumption that $a_i$ are all real comes in. (See \cite[Theorem
10.2.1]{Hille}.) Here, because we actually starts from an automorphic
function first, we do not need this assumption. For convenience of the
reader, we provide a complete proof of the proposition here.

\begin{proof} We first consider the analytic behavior of
$$
   Q(t)=-\frac1{2t'(\tau)^2}\left(\frac{t'''(\tau)}{t'(\tau)}
   -\frac32\left(\frac{t''(\tau)}{t'(\tau)}\right)^2\right)
$$
at a point $t_0=t(\tau_0)\neq\infty$ that does not correspond to any
elliptic point. We have
$$
  t(\tau)=t_0+c(\tau-\tau_0)+\cdots.
$$
for some constant $c$. Since $t(\tau)-t_0$ is also a
Hauptmodul, the constant $c$ cannot be $0$. From this we easily
see that $Q(t)$ is holomorphic at $t_0$.

We next consider the case $t_0=a_j\neq\infty$ corresponding
to an elliptic point $\tau_j$ of order $e_j$.
\begin{equation*}
\begin{split}
   t(\tau)=a_j+c_1(\tau-\tau_0)^{e_j}+c_2(\tau-\tau_0)^{e_j+1}+\cdots
\end{split}
\end{equation*}
for some constants $c_1,c_2\in\C$.

Again, because $t(\tau)-a_j$ is also a Hauptmodul, $c_1$ cannot be
equal to $0$. When $e_j\ge 3$, we have
\begin{equation*}
\begin{split}
  \frac{t'''(\tau)}{t'(\tau)}-\frac32\frac{t''(\tau)^2}{t'(\tau)^2}
  &=\frac{c_1e_j(e_j-1)(e_j-2)(\tau-\tau_0)^{e_j-3}+\cdots}{c_1e_j(\tau-\tau_0)^{e_j-1}+\cdots} \\
     &\qquad\qquad-\frac32\frac{c_1^2e_j^2(e_j-1)^2(\tau-\tau_0)^{2e_j-4}+\cdots}
       {c_1^2e_j^2(\tau-\tau_0)^{2e_j-2}+\cdots} \\
  &=\frac12\frac{1-e_j^2}{(\tau-\tau_0)^2}+\cdots.
\end{split}
\end{equation*}
When $e_j=2$, we have
\begin{equation*}
\begin{split}
  \frac{t'''(\tau)}{t'(\tau)}-\frac32\frac{t''(\tau)^2}{t'(\tau)^2}
  &=\frac{6c_2+\cdots}{2c_1(\tau-\tau_0)+\cdots}-\frac32
        \frac{4c_1^2+\cdots}{4c_1^2(\tau-\tau_0)^2+\cdots} \\
  &=\frac12\frac{1-e_j^2}{(\tau-\tau_0)^2}+\cdots.
\end{split}
\end{equation*}
Either way, we have
$$
  Q(t)=-\frac1{2c_1^2e_j^2(\tau-\tau_0)^{2e_j-2}+\cdots}\left(
   \frac12\frac{1-e_j^2}{(\tau-\tau_0)^2}+\cdots\right)
  =\frac{1-1/e_j^2}{4c_1^2(\tau-\tau_0)^{2e_j}}+\cdots,
$$
which implies that
$$
   Q(t)-\frac14\frac{1-1/e_j^2}{(t-a_j)^2},
$$
as a function of $t$, has at most a simple pole at $a_j$. Thus, if we
let $B_j$ denote the residue of $Q(t)$ at $a_j$, then
\begin{equation} \label{equation: P}
  P(t):=Q(t)-\frac14\sum_{j=1,a_j\neq\infty}^r\frac{1-1/e_j^2}{(t-a_j)^2}
  -\sum_{j=1,a_j\neq\infty}\frac{B_j}{t-a_j}
\end{equation}
will be a polynomial function in $t$. We now consider the behavior of
$Q(t)$ at $\infty$.

If $a_j\neq\infty$ for all elliptic points $\tau_j$, letting
$\tau_0\in\H$ be a point with $t(\tau_0)=\infty$, we have
$$
  t(\tau)=\frac c{\tau-\tau_0}+\cdots
$$
for some nonzero constant $c$. Then
\begin{equation} \label{equation: Q at infinity}
\begin{split}
  Q(t)&=-\frac1{2(-c(\tau-\tau_0)^{-2}+O(1))^2} \\
  &\qquad\qquad\times
    \left(\frac{-6c(\tau-\tau_0)^{-4}+O(1)}{-c(\tau-\tau_0)^{-2}+O(1)}
   -\frac32\frac{4c^2(\tau-\tau_0)^{-6}+O((\tau-\tau_0)^{-3})}
   {c^2(\tau-\tau_0)^{-4}+O((\tau-\tau_0)^{-2})}\right)
   \\
  &=O\left((\tau-\tau_0)^4\right)=O(t^{-4})
\end{split}
\end{equation}
as $\tau\to\tau_0$. In particular, we have $Q(t)\to 0$ as
$t\to\infty$, which implies that $P(t)=0$. Moreover, the coefficients
of $t^{-1}$, $t^{-2}$, $t^{-3}$ in the expansion of
\begin{equation} \label{equation: P at infinity}
  \frac14\sum_{j=1,a_j\neq\infty}^r\frac{1-1/e_j^2}{(t-a_j)^2}
 +\sum_{j=1,a_j\neq\infty}\frac{B_j}{t-a_j}
\end{equation}
at $t=\infty$ are
$$
  \sum_{j=1}^r B_j, \qquad
  \sum_{j=1}^r\left(a_jB_j+\frac14(1-1/e_j^2)\right), \qquad
  \sum_{j=1}^r\left(a_j^2B_j+\frac12a_j(1-1/e_j^2)\right),
$$
respectively. In view of \eqref{equation: Q at infinity}, all the
three sums must be equal to $0$. This proves the proposition in the
case $a_j\neq\infty$ for all $j$.

Likewise, if $a_j=\infty$ for some elliptic point $\tau_j$, say,
$a_r=\infty$, we have
$$
  t(\tau)=\frac c{(\tau-\tau_r)^{e_r}}+\cdots
$$
for some nonzero constant $c$. Then
\begin{equation*}
\begin{split}
  Q(t)&=-\frac1{2(-ce_r(\tau-\tau_r)^{-e_j-1}+\cdots)^2}
   \Bigg(\frac{-ce_r(e_r+1)(e_r+2)(\tau-\tau_0)^{-e_r-3}+\cdots}
   {-ce_r(\tau-\tau_r)^{-e_r-1}+\cdots} \\
  &\qquad\qquad\qquad-\frac32
   \frac{c^2e_r^2(e_r+1)^2(\tau-\tau_r)^{-2e_r-4}+\cdots}
   {c^2e_r^2(\tau-\tau_r)^{-2e_r-2}+\cdots}\Bigg) \\
  &=\frac{1-1/e_r^2}{4c^2}(\tau-\tau_r)^{2e_r}+\cdots
   =\frac{1-1/e_r^2}{4t^2}+O(t^{-3})
\end{split}
\end{equation*}
as $t\to\infty$. This also shows that the polynomial $P(t)$ in
\eqref{equation: P} is actually $0$. Also, comparing the coefficients
of $t^{-1}$ and $t^{-2}$ in the expansion of \eqref{equation: P at
  infinity} with the asymptotic of $Q(t)$ at $\infty$ given above, we
find
$$
  \sum_{j=1}^{r-1}B_j=0, \qquad
  \sum_{j=1}^{r-1}\left(a_jB_j+\frac14(1-1/e_j^2)\right)=\frac14(1-1/e_j^2).
$$
This proves the proposition for the case $a_j=\infty$ for some $j$.
\end{proof}

Here we will give two examples of automorphic derivatives. Let us
first fix some notations.

\begin{Notation} For an Eichler order $\O$ of level $N$ in an
  indefinite quaternion algebra of discriminant $D$ over $\Q$, we let
  $\Gamma_D(N)$ denote the group $\Gamma(\O)$ and $X_D(N)$ denote the
  Shimura curve associated to $\O$. If we let $W_D$ be the group of
  Atkin-Lehner involutions $\{w_e:e|D\}$, then $\Gamma^\ast_D(N)$ and
  $X^\ast_D(N)$ will denote the group $\bigcup_{e|D}w_e\Gamma_D(N)$
  and the quotient curve $X_D(N)/W_D$, respectively.

  Also, as is customary, if a Shimura curve has genus $g$ with $m_i$
  elliptic points of order $e_i$, we use the notation
  $(g;e_1^{m_1},\ldots,e_r^{m_r})$ to encode the signature
  of the curve.
\end{Notation}

\begin{Example} \label{example: X61}
  Consider the Shimura curve $X^\ast_6(1)$. The
  signature of $X^\ast_6(1)$ is $(0;2,4,6)$. Choose a
  Hauptmodul $t$ of $X^\ast_6(1)$ by requiring that $t$ takes values
  $0$, $1$, and $\infty$ at the elliptic points of orders $6$, $2$,
  and $4$, respectively. By Proposition \ref{proposition: Q},
  $$
    Q(t)=\frac{35/144}{t^2}+\frac{3/16}{(t-1)^2}+
    \frac{B_1}t+\frac{B_2}{t-1},
  $$
  where $B_1$ and $B_2$ satisfy
  $$
    B_1+B_2=0, \qquad B_2+\frac{35}{144}+\frac3{16}=\frac{15}{64}.
  $$
  Thus,
  $$
    Q(t)=\frac{35/144}{t^2}+\frac{3/16}{(t-1)^2}
    +\frac{113/576}t+\frac{113/576}{1-t}.
  $$
  The local exponents of the differential equation $d^2G/dt^2+Q(t)G=0$
  at $0$, $1$, and $\infty$, are $\{5/12,7/12\}$, $\{1/4,3/4\}$, and
  $\{-5/8,-3/8\}$, respectively. Thus, the differential equation
  satisfied by $t^{-5/12}(1-t)^{-1/4}t'(\tau)^{1/2}$ (as a function of
  $t$) will have local exponents $\{0,1/6\}$, $\{0,1/2\}$, and
  $\{1/24,7/24\}$ at $0$, $1$, and $\infty$, respectively. In other
  words, $t^{-5/12}(1-t)^{-1/4}t'(\tau)^{1/2}$, as a function of $t$,
  is a solution of the hypergeometric differential equation
  $$
    \theta\left(\theta-\frac16\right)F-t\left(\theta+\frac1{24}\right)
    \left(\theta+\frac7{24}\right)F=0, \quad\theta=t\frac d{dt}.
  $$
  Therefore, if we fix a representative $\tau_1\in\H$ of the elliptic
  point of $6$, then in a neighborhood of $\tau_1$, we have
  \begin{equation} \label{equation: t' for X61}
    t'(\tau)=t^{5/6}(1-t)^{1/2}\left(C_1{}_2F_1\left(
    \frac1{24},\frac7{24};\frac56;t\right)
    +C_2t^{1/6}{}_2F_1\left(\frac5{24},\frac{11}{24};\frac76;t\right)
    \right)^2,
  \end{equation}
  where $C_1$ and $C_2$ are complex numbers depending on the choice of
  the embedding of the quaternion algebra into $M(2,\R)$ and the choice
  of $\tau_1$. Note that from \eqref{equation: whatever}, we know that
  $C_1$ is nonzero.
  Then by Theorem \ref{theorem: basis}, for even positive integers
  $k\ge 4$, the automorphic forms whose $t$-expansions near the point
  $\tau_1$ are
  \begin{equation} \label{equation: basis for X61}
    t^{\{5k/12\}}(1-t)^{\{k/4\}}t^j \left({}_2F_1\left(
    \frac1{24},\frac7{24};\frac56;t\right)
    +Ct^{1/6}{}_2F_1\left(\frac5{24},\frac{11}{24};\frac76;t\right)
    \right)^k,
  \end{equation}
  $j=0,\ldots,\lfloor 5k/12\rfloor+\lfloor k/4\rfloor+\lfloor
  3k/8\rfloor-k$, form a basis for the space of automorphic forms of
  weight $k$ on $\Gamma^\ast_6(1)$, where $C=C_2/C_1$ and for a
  rational number $x$ we let $\{x\}=x-\lfloor x\rfloor$ denotes the
  fractional part of $x$. In Section \ref{section: 246}, we will
  compute Hecke operators relative to this basis.
\end{Example}

More generally, automorphic forms on triangle groups can all be
expressed in terms of hypergeometric functions.

\begin{Theorem} \label{theorem: triangle}
  Assume that a Shimura curve $X$ has signature
  $(0;e_1,e_2,e_3)$. Let $t(\tau)$ be the Hauptmodul of $X$ with
  values $0$, $1$, and $\infty$ at the elliptic points of order $e_1$,
  $e_2$, and $e_3$, respectively. Let $k\ge 4$ be an even integer.
  Then a basis for the space of automorphic forms of weight $k$ on $X$
  is given by
  $$
    t^{\{k(1-1/e_1)/2\}}(1-t)^{\{k(1-1/e_2)/2\}}t^j\left(
    {}_2F_1(a,b;c;t)+Ct^{1/e_1}{}_2F_1(a',b',c';t)\right)^k,
  $$
  $j=0,\ldots,\lfloor k(1-1/e_1)/2\rfloor+\lfloor k(1-1/e_2)/2\rfloor
  +\lfloor k(1-1/e_3)/2\rfloor-k$, for some constant $C$, where
  $$
    a=\frac12\left(1-\frac1{e_1}-\frac1{e_2}-\frac1{e_3}\right), \qquad
    b=a+\frac1{e_3}, \qquad c=1-\frac1{e_1}
  $$
  and
  $$
    a'=a+\frac1{e_1}, \qquad b'=b+\frac1{e_1}, \qquad
    c'=c+\frac2{e_1}.
  $$
\end{Theorem}

\begin{proof} The proof follows the same argument as in the preceding
  example. Here we just remind the reader that the hypergeometric
  function $_2F_1(a,b,c;t)$ is a solution of
  $$
    \theta(\theta+c-1)F-t(\theta+a)(\theta+b)F=0, \qquad
    \theta=\frac d{dt},
  $$
  whose local exponents at $0$, $1$, and $\infty$ are $\{0,1-c\}$,
  $\{0,c-a-b\}$, and $\{a,b\}$, respectively. From this, it is easy to
  figure out the parameters in the hypergeometric functions. We omit
  the details.
\end{proof}

\begin{Example}
  Consider the Shimura curve $X^\ast_{10}(1)$. The
  signature of $X^\ast_{10}(1)$ is $(0;3,2^3)$. In \cite{Elkies}, Elkies
  showed that there is a Hauptmodul on $X^\ast_{10}(1)$ with values $0$
  at the elliptic point of order $3$ and values $2$, $27$, and
  $\infty$ at the three elliptic points of order $2$. Moreover, using
  the covering $X^\ast_{10}(3)\to X^\ast_{10}(1)$, he also showed that an
  automorphic differential equation for $X^\ast_{10}(1)$ is
  $$
    t(t-2)(t-27)F''+\frac{10t^2-203t+216}6F'+\frac{7t-56}{144}F=0,
  $$
  which is the same as
  $$
    \theta^2F+\frac{4t^2-29t-108}{6(t-2)(t-27)}\theta F
   +\frac{7t(t-8)}{144(t-2)(t-27)}F=0, \quad \theta=t\frac d{dt}.
  $$
  Using Proposition \ref{proposition: Schwarzian}, we find that the
  automorphic derivative $D(t,\tau)$ associated to $t$ is
  \begin{equation} \label{equation: Q10}
    D(t,\tau)=\frac{3t^4-119t^3+3157t^2-7296t+10368}{16t^2(t-2)^2(t-27)^2}.
  \end{equation}
  Here we will use Propositions \ref{proposition: D} and
  \ref{proposition: Q} to obtain the same result.

  According to Proposition \ref{proposition: Q}, we have
  $$
    D(t,\tau)=\frac{2/9}{t^2}+\frac{3/16}{(t-2)^2}+\frac{3/16}{(t-27)^2}
    +\frac{B_1}t+\frac{B_2}{t-2}+\frac{B_3}{t-27},
  $$
  where $B_i$ satisfy
  $$
    B_1+B_2+B_3=0, \qquad 2B_2+27B_3+\frac29+\frac3{16}
   +\frac3{16}=\frac3{16},
  $$
  which yield
  $$
    B_1=\frac{25}2B_3+\frac{59}{288}, \qquad
    B_2=-\frac{27}2B_3-\frac{59}{288}.
  $$

  On the other hand, the Shimura curve $X_{10}^\ast(3)$ has signature
  $(0;2^4,3)$. According to \cite{Elkies}, there is a Hauptmodul $u$
  on
  $X_{10}^\ast(3)$ whose relation with $t$ is
  \begin{equation} \label{equation: 10 t-u}
    t=\frac{216(u-1)^3}{(u+1)^2(9u^2-10u+17)}
  \end{equation}
  and whose values at the four elliptic points of order $2$ and the
  elliptic point of order $3$ are $(5\pm 2\sqrt{-5})/9$,
  $(5\pm 8\sqrt{-2})/9$, and $\infty$, respectively. (Note that in
  (57) of \cite{Elkies}, the factor $9x^2-10x+17$ in the denominator was
  misprinted as $9x^2-10x+7$.) Also, the action of
  Atkin-Lehner involution $w_3$ is
  \begin{equation} \label{equation: 10 w3}
    w_3:u\longmapsto \frac{10}9-u.
  \end{equation}
  By Proposition \ref{proposition: Q},
  \begin{equation} \label{equation: Du}
  \begin{split}
    D(u,\tau)&=\frac{3(162u^2-180u+10)/16}{(9u^2-10u+5)^2}
    +\frac{3(162u^2-180u-206)/16}{(9u^2-10u+17)^2} \\
    &\qquad\qquad+\frac{C_1+C_2u}{9u^2-10u+5}
    +\frac{C_3+C_4u}{9u^2-10u+17}
  \end{split}
  \end{equation}
  for some constants $C_1,\ldots,C_4$. According to the proof of
  Proposition \ref{proposition: Q}. The constants $C_1,\ldots,C_4$
  satisfy $Q(u)=2/9t^{2}+O(t^{-3})$ as $t\to\infty$. Thus,
  $$
     C_4=-C_2, \qquad C_3=-C_1-\frac{19}4.
  $$
  Now observe that the Atkin-Lehner involution $w_3:u\mapsto10/9-u$
  switches the two roots of $9u^2-10u+5$ and the two roots of
  $9u^2-10u+17$ and fixes $\infty$. From this we deduce that
  $D(u,\tau)=D(10/9-u,\tau)$, which implies that
  the right-hand side of \eqref{equation: Du} is invariant under the
  substitution $u\mapsto10/9-u$. From this, we infer that $C_2=C_4=0$.
  That is,
  \begin{equation*}
  \begin{split}
    D(u,\tau)&=\frac{3(162u^2-180u+10)/16}{(9u^2-10u+5)^2}
    +\frac{3(162u^2-180u-206)/16}{(9u^2-10u+17)^2} \\
    &\qquad\qquad
    +\frac{C_1}{9u^2-10u+5}+\frac{-C_1-19/4}{9u^2-10u+17}.
  \end{split}
  \end{equation*}
  Now let $R(x)=216(x-1)^3/(x+1)^2/(9x^2-10x+17)$. We have $t=R\circ
  u$. Thus, by Part (2) of Proposition \ref{proposition: D}, we have
  $$
    D(t,\tau)=D(R(u),u)+\frac{D(u,\tau)}{(dR(u)/du)^2}.
  $$
  Expressing the left-hand side in terms of $u$ and comparing it with
  the right-hand side, we find that the constants $B_3$ and $C_1$
  are
  $$
    B_3=-\frac{953}{97200}, \qquad C_1=-\frac1{24}.
  $$
  This gives us \eqref{equation: Q10}.
\end{Example}

\begin{Remark} In \cite{Tu}, Tu determined the Schwarzian differential
  equations associated to $X_D^\ast(N)$ for the following pairs of
  $(D,N)$:
  \begin{equation*}
  \begin{split}
   &(6,1),\ (6,5), \ (6,7), \ (6,13), \ (10,1), \ (10,3), \ (10,7), \\
   &(14,1), \ (14,3), \ (14,5), \ (15,1), \ (15,2), \ (15,4), \
   (21,1),\\
   &(21,2), \ (26,1), \ (26,3), \ (35,1), \ (35,2), \ (39,1), \ (39,2).
  \end{split}
  \end{equation*}
  For example, for $(D,N)=(39,1)$, she showed that there exists a
  Hauptmodul $t(\tau)$ for $X_{39}^\ast(1)$ that takes values $\pm
  2i$, $(-1\pm\sqrt{-3})/2$, and $(-23+\sqrt{-3})/14$ at the CM-points
  of discriminants $-52$, $-39$, and $-156$, respectively, and the
  function $Q(t)$ in the Schwarzian differential equation associated
  to $t(\tau)$ is
  $$
    Q(t)=-\frac{3(2596+7104t+9692t^2+12348t^3+13149t^4+9522t^5
    +4367t^6+1086t^7+97t^8)}
    {4(4+t^2)^2(1+t+t^2)^2(19+23t+7t^2)^2}.
  $$
\end{Remark}

The rest of the paper will be devoted to the computation of Hecke
operators on automorphic forms. But before we work on the case of
Shimura curves, let us first work out a familiar example from classical
modular curves to give the reader a clearer idea about our
approach.
\end{section}

\begin{section}{Computing Hecke operators -- an example}
Let $\Delta(\tau)=\eta(\tau)^{24}$ be the unique normalized Hecke
eigenform on $\SL(2,\Z)$. Of course, from the Fourier expansion
$\Delta(\tau)=q-24q^2+\cdots$, we immediately see that the eigenvalue
for the Hecke operator $T_2$ is $-24$. Our goal in this section
is to obtain the same result without resorting to Fourier expansions.

It is a classical identity that
$$
  \Delta(\tau)=t{}_2F_1\left(\frac1{12},\frac5{12};1;1728t\right)^{12},
  \qquad t(\tau)=\frac1{j(\tau)},
$$
where $j(\tau)$ is the elliptic $j$-function. (In fact, this can also
be verified using a slightly modified version of our Theorem
\ref{theorem: triangle}.) Thus, as long as the imaginary part of
$\tau$ is large, we may expand $\Delta(\tau)$ with respect to $t$.
Assuming that $\Im\tau$ is large, by the definition of $T_2$, we have
\begin{equation} \label{equation: T2Delta}
\begin{split}
  T_2\Delta(\tau)&=2^{11}t(2\tau)
  {}_2F_1\left(\frac1{12},\frac5{12};1;t(2\tau)\right)^{12}
 +\frac12t(\tau/2)
  {}_2F_1\left(\frac1{12},\frac5{12};1;t(\tau/2)\right)^{12}\\
  &\qquad\qquad
 +\frac12t((\tau+1)/2)
  {}_2F_1\left(\frac1{12},\frac5{12};1;t((\tau+1)/2)\right)^{12}.
\end{split}
\end{equation}
Now suppose that we are allowed to use Fourier expansions for the
moment. We have $t(\tau)=1/j(\tau)=q-744q^2+\cdots$ and
\begin{equation} \label{equation: T2Delta 1}
\begin{split}
  t(2\tau)&=q^2-744q^4+\cdots=t(\tau)^2+1488t(\tau)^3+\cdots, \\
  t(\tau/2)&=q^{1/2}-744q+\cdots=t(\tau)^{1/2}-744t(\tau)+\cdots, \\
  t((\tau+1)/2)&=-q^{1/2}-744q+\cdots=-t(\tau)^{1/2}-744t(\tau)+\cdots.
\end{split}
\end{equation}
Substituting these expressions into \eqref{equation: T2Delta}, we get
\begin{equation*}
\begin{split}
  T_2\Delta(\tau)&=2^{11}\left(t^2+\cdots\right)
     \left(1+60t^2+\cdots\right)^{12} \\
  &\qquad\qquad+\frac12\left(t^{1/2}-744t+\cdots\right)
     \left(1+60t^{1/2}+\cdots\right)^{12} \\
  &\qquad\qquad+\frac12\left(-t^{1/2}-744t+\cdots\right)
     \left(1-60t^{1/2}+\cdots\right)^{12} \\
  &=-24t(\tau)+\cdots=-24\Delta(\tau).
\end{split}
\end{equation*}
From this, we see that the eigenvalue for $T_2$ is indeed $-24$.
Note that there is an ambiguity in the choice of the square root of $t$
in \eqref{equation: T2Delta 1}, but it does not affect the final
result.

Of course, we have cheated a little bit in the above computation by
using Fourier expansion in \eqref{equation: T2Delta 1}. We now discuss
how to obtain the same $t$-expansions without using $q$-expansions.
The idea is to use the so-called \emph{modular equation}, which is the
polynomial relation satisfied by $j(\tau)$ and $j(2\tau)$.

Observe that $t(\tau)$ and $t(2\tau)$ are both modular functions on
$\Gamma_0(2)$. Let $u(\tau)$ be a Hauptmodul of $\Gamma_0(2)$. Since
$X_0(2)\to X_0(1)$ is a covering of degree $3$, we have
$t(\tau)=R(u(\tau))$ for some rational function $R$ of exactly degree
$3$. Now $t(2\tau)=t(-1/2\tau)=R(u(-1/2\tau))$. Since $\SM0{-1}20$
normalizes $\Gamma_0(2)$, $u(-1/2\tau)$ is also a Hauptmodul and
therefore
$$
  u(-1/2\tau)=\frac{au(\tau)+b}{cu(\tau)+d}
$$
for some $a,b,c,d\in\mathrm{GL}(2,\C)$. Hence,
$$
  t(2\tau)=R\left(\frac{au(\tau)+b}{cu(\tau)+d}\right).
$$
In other words, the polynomial relation between $t(\tau)$ and
$t(2\tau)$ is just the relation between $t=R(u)$ and
$s=R((au+b)/(cu+d))$.

Now $t$ has values $0$, $1/1728$, and $\infty$ at the cusp $P_\infty$,
the elliptic point $P_2$ of order $2$ and the elliptic point $P_3$ of
order $3$, respectively. Above these three points, we have the
following ramification data

\centerline{
\begin{tikzpicture}[scale=0.07]
  \node[place] (Pi) [label=below:$P_\infty$, token=1] at (-30,0) {};
  \node[place] (P2) [label=below:$P_2$, token=1] at (10,0) {};
  \node[place] (P3) [label=below:$P_3$, token=1] at (40,0) {};
  \node[place] (Qi1) [label=above:$Q_\infty$, token=1] at (-40,20) {};
  \node[place] (Qi2) [label=above:$Q'_\infty$, token=1] at (-20,20)
  {};
  \node[place] (Q21) [label=above:$Q_2$, token=1] at (0,20) {};
  \node[place] (Q22) [label=above:$Q'_2$, token=1] at (20,20) {};
  \node[place] (Q3) [label=above:$Q_3$, token=1] at (40,20) {};
  \draw[line width=.6pt] (Qi1) --  node[left=2pt] {$1$} (Pi);
  \draw[line width=.6pt] (Qi2) --  node[right=2pt] {$2$}(Pi);
  \draw[line width=.6pt] (Q21) --  node[left=2pt] {$1$}(P2);
  \draw[line width=.6pt] (Q22) --  node[right=2pt] {$2$}(P2);
  \draw[line width=.6pt] (Q3) --  node[right=2pt] {$3$}(P3);
\end{tikzpicture}
}
\noindent{Here the numbers next to the lines are the ramification
  indices.}

Choose the Hauptmodul $u$ of $\Gamma_0(2)$ with values $u(Q_\infty)=0$,
$u(Q_2)=1$, and $u(Q_3)=\infty$. From the ramification data at
$P_\infty$ and $P_3$, we have $R(u)=Au(u-\alpha)^2$ for some
$\alpha\in\C$. Also, the ramification data $P_2$ implies
$Au(u-\alpha)^2-1/1728=A(u-1)(u-\beta)^2$ for some $\beta\in\C$.
Comparing the coefficients we find $A=1/108$, $\alpha=3/4$, and
$\beta=1/4$. Furthermore, the Atkin-Lehner involution $w_2$ switches
the two cusps $Q_\infty$ and $Q_\infty'$ and fixes the elliptic point
$Q_2$ of order $2$. Thus,
$$
  w_2:u\longmapsto\frac{4u-3}{5u-4}.
$$
Eliminating $u$, we find that the relation between
$t=R(u)=u(u-3/4)^2/108$ and $s=R((4u-3)/(5u-4))$ is
\begin{equation*}
\begin{split}
  \Phi_2(s,t)&=s^3+t^3-st+1488s^2t-162000s^3t+1488st^2 \\
  &\qquad+40773375s^2t^2+8748000000s^3t^2-162000st^3 \\
  &\qquad+8748000000s^2t^3-157464000000000s^3t^3.
\end{split}
\end{equation*}
Solving $\Phi_2(s,t)=0$ for $s$, we find the three roots are
\begin{equation*}
\begin{split}
  s&=t^2+1488t^3+2053632t^4+\cdots, \\
  s&=t^{1/2}-744t+357024t^{3/2}+\cdots, \\
  s&=-t^{1/2}-744t-357024t^{3/2}+\cdots,
\end{split}
\end{equation*}
which agree with the $t$-expansions of $t(2\tau)$, $t(\tau/2)$, and
$t((\tau+1)/2)$ given in \eqref{equation: T2Delta 1}.

Indeed, using differential equations and modular equations, we can
compute Hecke operators on the spaces of modular forms on $\SL(2,\Z)$
without resorting to Fourier expansions. In the next two sections, we
will use the same idea to compute Hecke operators in the case of
Shimura curves.
\end{section}

\begin{section}{Hecke operators on $X_6^\ast(1)$}
  \label{section: 246}

Hecke operators on the space of automorphic forms on Shimura curves
associated to Eichler orders are defined in the same way as in the
case of classical modular curves. For simplicity, we assume that the
quaternion algebra $B$ is over $\Q$ and has discriminant $D$. Fix an
embedding $\iota: B\to M(2,\R)$. Let $\O$ be an Eichler order of level
$N$. For a prime $p$ not dividing $DN$, we pick an
element of reduced norm $p$ in $\O$. Then one can show that
$\Gamma(\O)$ and $\iota(\alpha)^{-1}\Gamma(\O)\iota(\alpha)$ are
commensurable so that
$\Gamma(\O)\backslash\Gamma(\O)\iota(\alpha)\Gamma(\O)$ has finitely
many right cosets. (In fact, the number of right cosets is
$p+1$.) Then for an automorphic form $f(\tau)$ of even weight $k$ on
$\Gamma(\O)$, the action of Hecke operator $T_p$ on $f(\tau)$ is
defined by
\begin{equation} \label{equation: Tp}
  T_p:f\longmapsto p^{k/2-1}\sum_{\gamma\in
  \Gamma(\O)\backslash\Gamma(\O)\iota(\alpha)\Gamma(\O)}
  \frac{(\det\gamma)^{k/2}}{(c\tau+d)^k}f(\gamma\tau),
\end{equation}
where for a coset representative $\gamma$, we write $\gamma=\SM
abcd$. Hecke operators $T_n$ for general $n$ with $(n,DN)=1$ are
slightly more complicated.

In this section, we will compute Hecke operators on automorphic forms
on $X_6^\ast(1)$. The computation follows that in the previous section
in principle, but several issues arise.
\begin{enumerate}
\item Proposition \ref{proposition: Schwarzian} only says that
  $t'(\tau)^{1/2}$ satisfies the Schwarzian differential equation, but
  it does not say which solution corresponds to $t'(\tau)^{1/2}$. This
  is not a problem in the example in the previous section because the
  hypergeometric differential equation
  $\theta^2F-1728t(\theta+1/12)(\theta+5/12)F=0$ has a unique
  solution (up to scalars) that is holomorphic at $t=0$ and
  $t^{-1/2}t'(\tau)^{1/2}$ must be a multiple of this solution. (The other
  solutions have a logarithmic singularity at $t=0$.) Here we need to
  find two linearly independent solutions of the differential equation
  and then find an appropriate linear combination that corresponds to
  $t'(\tau)^{1/2}$.
\item Unlike the example in the previous section, here we also need to
  find the $t$-expansion for $\tau$. Nonetheless, this problem is
  relatively simple to settle once the first problem is answered.
\item In the example in the previous section, since the $t$-expansions
  converge only for $\tau$ with large imaginary parts, there is an
  obvious choice of coset representatives $\SM2001$, $\SM1002$, and
  $\SM1102$, but here it is not immediately clear how we should choose
  coset representatives.
\item Even if we are able to find the polynomial relation
  $\Phi(s,t)=0$ between $t(\tau)$ and $s(\tau)=t(\gamma\tau)$,
  $\gamma\in\Gamma(\O)\iota(\alpha)\Gamma(\O)$, and solve the equation
  for $s$ as $t$-series, we still need to determine which solution of
  the equation is matched with which coset representatives. In the
  example in the previous section, this is relatively simple. The
  solution starting with $t^2+\cdots$ must correspond to $\SM2001$,
  while it does not really matter how the other two solutions are
  matched with coset representatives.
\item Unlike the example in the previous section where the choice of
  coset representatives makes $t(\gamma\tau)\to 0$ as $t(\tau)\to 0$
  so that to find the $t$-expansion of $f(\gamma\tau)$, we only need
  to substitute $t$ by the $t$-expansion of $t(\gamma\tau)$ in $f$,
  here we also need to find a method to determine the $t$-expansion of
  $f(\gamma\tau)$ for each coset representative. This is perhaps the
  most complicated part of the computation. The Jacquet-Langlands
  correspondence will be very useful in this part. We now recall an
  explicit version of the correspondence in the case of quaternion
  algebras over $\Q$.
\end{enumerate}

\begin{Proposition}[{\cite{Jacquet-Langlands,Shimizu}}]
  \label{proposition: Jacquet-Langlands}
  Let $D$ be discriminant of an indefinite quaternion algebra over
  $\Q$. Let $N$ be a positive integer relatively prime to $D$. For an
  Eichler order $\O(D,N)$ of level $(D,N)$ and a positive even
  integer, let $S_k(\O(D,N))$ denote the space of automorphic forms on
  $\Gamma(\O(D,N))$. Then
  $$
    S_k(\O(D,N))\simeq
    S_k^{D\text{-\rm{new}}}(DN):=\bigoplus_{d|N}\bigoplus_{m|N/d}
    S_k^\new(dD)^{[m]}
  $$
  as Hecke modules. Here
  $$
    S_k^\new(dD)^{[m]}=\{f(m\tau):~f(\tau)\in S_k^\new(dD)\}
  $$
  and $S_k^\new(dD)$ denotes the newform subspace of
  cusp forms of weight $k$ on $\Gamma_0(dD)$. In other words, for each
  Hecke eigenform $f(\tau)$ in $S_k^{D\text{-\rm{new}}}(DN)$, there
  corresponds a Hecke eigenform $\widetilde f(\tau)$ in
  $S_k(\O(D,N))$ that shares the same Hecke eigenvalues. Moreover, for
  a prime divisor $p$ of $D$, if the Atkin-Lehner involution $W_p$
  acts on $f$ by $W_pf=\epsilon_p f$, then
  $$
    W_p\widetilde f=-\epsilon_p\widetilde f.
  $$
\end{Proposition}

We now consider the case $X_6^\ast(1)$. According to Example
\ref{example: X61}, if we choose the Hauptmodul $t$ with values $0$,
$1$, and $\infty$ at the elliptic points of order $6$, $2$, and $4$,
respectively, then the space of automorphic forms of weight $k$ has a
basis given by \eqref{equation: basis for X61}. Here we rescale the
Hauptmodul $t$ such that it has values $0$, $-540$, and $\infty$ at
the elliptic points of order $6$, $2$, and $4$, respectively. (The
purpose of this scaling is to make the coefficients of $t$-series in
the future computation simpler.) Then the
basis for the space of automorphic forms become
\begin{equation} \label{equation: basis for X61 2}
\begin{split}
 &g_\ell=t^{\{5k/12\}}(1+t/540)^{\{k/4\}}t^\ell \\
 &\qquad\times\left(
  {}_2F_1\left(\frac1{24},\frac7{24};\frac56;-\frac t{540}\right)
 -Ct^{1/6} {}_2F_1\left(\frac5{24},\frac1{24};\frac76;-\frac t{540}
  \right)\right)^k,
\end{split}
\end{equation}
$\ell=0,\ldots,\lfloor 5k/12\rfloor+\lfloor k/4\rfloor
+\lfloor 3k/8\rfloor-k$, where $C$ is a nonzero constant. We will
compute Hecke operators relative to this basis.

Let us first fix the quaternion algebra $B$ of
discriminant $6$ to be $\JS{-1,3}\Q$, i.e., the algebra generated by
$I$ and $J$ over $\Q$ with the relations
$$
  I^2=-1, \qquad J^2=3, \qquad IJ=-JI,
$$
and choose the embedding
$\iota:B\to M(2,\R)$ to be
$$
  I\longmapsto\M0{-1}10, \qquad
  J\longmapsto\M{\sqrt3}00{-\sqrt3}
$$
as in \cite[Section 5.5.2]{Alsina-Bayer}. Fix the maximal order $\O$
to be $\Z+\Z I+\Z J+\Z(1+I+J+IJ)/2$. Then
$$
  \Gamma(\O)=\left\{\frac12\M\alpha\beta{-\beta'}{\alpha'}
  \in\SL(2,\R):~\alpha,\beta\in\Z[\sqrt3],~\alpha\equiv\beta\mod 2\right\},
$$
where $\alpha'$ and $\beta'$ denote the Galois conjugates of $\alpha$
and $\beta$, respectively.

As in \cite[Section 5.5.2]{Alsina-Bayer}, we choose the
representatives of elliptic points of order $2$, $4$, $6$ by 
$$
  P_2=(\sqrt 6-\sqrt 2)i/2, \qquad P_4=i, \qquad
  P_6=(-1+i)/(1+\sqrt3)
$$
with the isotropy subgroups generated by
$$
  M_2=\frac1{\sqrt6}\M0{-3+\sqrt3}{3+\sqrt3}0, \qquad
  M_4=\frac1{\sqrt2}\M1{-1}11,
$$
and
\begin{equation} \label{equation: A6}
  M_6=\frac1{2\sqrt3}\M{3+\sqrt3}{3-\sqrt3}{-3-\sqrt3}{3-\sqrt3},
\end{equation}
respectively. A fundamental domain for $X^\ast_6(1)$ is given by

\centerline{\epsfig{file=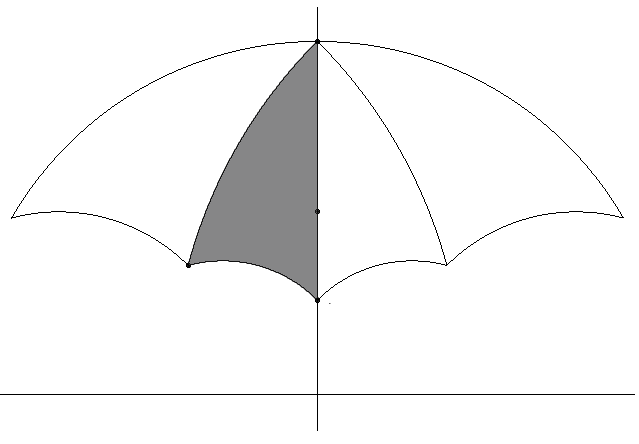,height=2.1in,width=3in}}

\noindent Here the grey area represents a fundamental domain for
$X_6^\ast(1)$. The four marked points on the boundary are
$$
  P_4=i, \quad P_6=\frac{-1+i}{1+\sqrt3}, \quad (2-\sqrt3)i, \quad
  P_2=\frac{(\sqrt6-\sqrt2)i}2,
$$
respectively. (Note that the action of $M_2$ maps $P_4$ to
$(2-\sqrt3)i$.) The grey area and three other white areas form a
fundamental domain for $X_6(1)$. (See \cite[Figure 5.1]{Alsina-Bayer}
and \cite{Voight-computation}.)

To compute $T_5$ on the functions in \eqref{equation: basis for X61},
we need to choose appropriate coset representatives $\gamma_j$,
$j=1,\ldots,6$, for $\Gamma^\ast(\O)\backslash\Gamma^\ast(\O)
\SM1{-2}21\Gamma^\ast(\O)$, where $\SM1{-2}21$ is the image of the
element $1+2I$ of reduced norm $5$ in $\O$ under $\iota$. At the
hindsight, if our goal is just to compute Hecke operators, it does not
really matter how we choose $\gamma_j$, as long as $\gamma_jP_6$ are
the same point for $j=1,\ldots,6$. Here we take the somehow natural
choice having the property that $\gamma_jP_6$ is in the fundamental
domain given above.

\begin{Lemma} \label{lemma: X61 cosets}
  Let the notations be given as above. A complete set of
  right coset representatives of $\Gamma^\ast(\O)\backslash
\Gamma^\ast(\O)\SM1{-2}21\Gamma^\ast(\O)$ is given by
\begin{equation*}
\begin{split}
  \gamma_0=\frac1{\sqrt{6}}\M{6+\sqrt3}{\sqrt3}{\sqrt3}{6-\sqrt3},
 &\qquad
  \gamma_1=\frac1{\sqrt{2}}\M{3+\sqrt3}{2}{-2}{3-\sqrt3}, \\
  \gamma_2=\frac1{\sqrt{6}}\M{3+2\sqrt3}{6-\sqrt3}{-6-\sqrt3}{3-2\sqrt3},
 &\qquad
  \gamma_3=\frac1{\sqrt{2}}\M{\sqrt3}{4-\sqrt3}{-4-\sqrt3}{-\sqrt3},
  \\
  \gamma_4=\frac1{\sqrt{6}}\M{-3+\sqrt3}{6-2\sqrt3}{-6-2\sqrt3}{-3-\sqrt3},
 &\qquad
  \gamma_5=\frac1{\sqrt{2}}\M{-3}{2-\sqrt3}{-2-\sqrt3}{-3}.
\end{split}
\end{equation*}
These coset representatives have the properties that
\begin{equation} \label{equation: gammaj}
  \gamma_j=\gamma_0M_6^j
\end{equation}
for all $j$ and
$$
  \gamma_jP_6=\frac{-1+5i}{1+3\sqrt3}
$$
is in the fundamental domain given above for all $j$, where $M_6$ is
given in \eqref{equation: A6}. (The indices are arranged such that
$\gamma_0P_2,\ldots,\gamma_5P_2$ are located counterclockwise around
$\gamma_jP_6$.) Moreover, letting $\gamma_j=\SM{a_j}{b_j}{c_j}{d_j}$
and $z=e^{2\pi i/24}$, we have
\begin{equation} \label{equation: X61 ctau+d}
  c_jP_6+d_j=z^{-2j}(z+z^5-z^7).
\end{equation}
\end{Lemma}

\begin{proof} Everything can be verified by a direct computation. We
  omit the details.
\end{proof}

We next determine the expansion of $\tau$ as a $t$-series in a neighborhood of
$P_6$. In the lemma below, the sixth root $t(\tau)^{1/6}$ of $t(\tau)$ is
defined in a neighborhood of $P_6$ such that it becomes a holomorphic
function of $\tau$ near $P_6$ and takes positive real values along the
boundary of the fundamental domain from $P_6$ to $P_4$. Note that
In view of $t(M_6\tau)=t(\tau)$, we have $t(M_6\tau)^{1/6}=\epsilon
t(\tau)^{1/6}$ for some sixth root of unity $\epsilon$. Since the
function $\tau\to t(\tau)$ preserves orientation and is locally
$6$-to-$1$ at $P_6$, this root of unity is actually $e^{2\pi i/6}$. In
other words, we have
\begin{equation} \label{equation: sixth root}
  t(M_6^j\tau)^{1/6}=e^{2\pi ij/6}t(\tau)^{1/6}.
\end{equation}
Similarly, the function $(1+t/540)^{1/2}$ is defined in a way such
that it becomes a holomorphic function near $P_2$ and takes positive
values along the boundary from $P_2$ to $P_6$ and from $P_6$ to $P_4$.
Note that we have
$$
  (1+t(M_2\tau)/540)^{1/2}=-(1+t(\tau)/540)^{1/2},
$$
even though this fact is not needed in the sequel.

\begin{Lemma} \label{lemma: X61 t'}
Let
$$
  F_1={}_2F_1\left(\frac1{24},\frac7{24};\frac56;-\frac t{540}\right),
  \qquad
  F_2=t^{1/6}{}_2F_1\left(\frac5{24},\frac{11}{24};\frac76;-\frac t{540}\right)
$$
be two linearly independently solution of
\begin{equation} \label{equation: DE for X61}
  \theta\left(\theta-\frac16\right)+\frac t{540}
    \left(\theta+\frac1{24}\right)\left(\theta+\frac7{24}\right)=0,
  \qquad \theta=t\frac d{dt}.
\end{equation}
We have
\begin{equation} \label{equation: C}
  \frac{\tau-P_6}{\tau-\overline P_6}=C\frac{F_2}{F_1}, \qquad
  C=\frac{P_2-P_6}{P_2-\overline P_6}\frac{e^{\pi i/6}}{\sqrt[6]{540}}
  \frac{\Gamma(5/6)\Gamma(17/24)\Gamma(23/24)}
  {\Gamma(7/6)\Gamma(13/24)\Gamma(19/24)}.
\end{equation}
Moreover, we have
\begin{equation} \label{equation: X61 t'}
\begin{split}
  t'(\tau)&=\frac{6t^{5/6}(1+t/540)^{1/2}}{C(P_6-\overline P_6)} \\
  &\qquad\times\left(
  {}_2F_1\left(\frac1{24},\frac7{24};\frac56;-\frac t{540}\right)
 -Ct^{1/6} {}_2F_1\left(\frac5{24},\frac1{24};\frac76;-\frac t{540}
  \right)\right)^2.
\end{split}
\end{equation}
That is, the constants $C$ in \eqref{equation: basis for X61 2} and
\eqref{equation: C} are the same.
\end{Lemma}

\begin{proof} The existence of a constant $C$ such that
  \eqref{equation: C} holds is well-known in the classical theory of
  automorphic functions. (Cf. Equations (48) and (49) of \cite{Elkies}.)
  Here we sketch a proof.

  From Example \ref{example: X61}, we know that
  $$
    t^{-5/12}(1+t/540)^{-1/4}t'(\tau)^{1/2}, \qquad
    \tau t^{-5/12}(1+t/540)^{-1/4}t'(\tau)^{1/2}
  $$
  are both solutions of the same differential equation
  \eqref{equation: DE for X61}. Thus,
  \begin{equation} \label{equation: tau 1}
    \tau=\frac{aF_1+bF_2}{cF_1+dF_2}
  \end{equation}
  for some complex numbers $a,b,c,d$. Now let $\gamma$ be a
  generator of the isotropy subgroup for 
  $P_6$. It is an elementary computation to show that
  \begin{equation} \label{equation: tau 2}
    \frac{\gamma\tau-P_6}{\gamma\tau-\overline P_6}=
    \epsilon\frac{\tau-P_6}{\tau-\overline P_6}
  \end{equation}
  for some primitive $6$th root of unity. The two facts \eqref{equation:
    tau 1} and \eqref{equation: tau 2} together imply that
  $(\tau-P_6)/(\tau-\overline{P}_6)=CF_1/F_2$ or
  $(\tau-P_6)/(\tau-\overline P_6)=CF_2/F_1$ for some nonzero complex
  number $C$. Since the left-hand side approaches $0$ as $\tau\to
  P_6$, it must be the second possibility that occurs. We then let
  $\tau\to P_2$ and use Gauss' formula
  $$
    {}_2F_1(a,b;c;1)=\frac{\Gamma(c)\Gamma(c-a-b)}{\Gamma(c-a)\Gamma(c-b)}
  $$
  to get the value of $C$. We now prove \eqref{equation: X61 t'}.

  From \eqref{equation: C}, we have
  $$
    \tau=\frac{P_6F_1-C\overline P_6F_2}{F_1-CF_2}.
  $$
  Differentiating the two sides with respect to $t$, we get
  $$
    \frac{d\tau}{dt}=C(P_6-\overline P_6)\frac{F_1dF_2/dt-F_2dF_1/dt}
    {(F_1-CF_2)^2}.
  $$
  Recall the formula that if $f_1$ and $f_2$ are two linearly
  independent solution of a second-order Fuchsian differential equation
  $d^2f/dt^2+p_1(t)df/dt+p_2(t)f=0$, then
  $$
    f_1\frac{df_2}{dt}-f_2\frac{df_1}{dt}=c\exp\left(-\int^t p_1(t)\,dt\right)
  $$
  for some constant $c$. Here we have
  $$
    p_1(t)=\frac{2(2t+675)}{3t(t+540)}
  $$
  and thus
  $$
    F_1\frac{dF_2}{dt}-F_2\frac{dF_1}{dt}=ct^{-5/6}(1+t/540)^{-1/2}
  $$
  for some $c$. Considering the leading coefficients, we find
  $c=1/6$. From this, we get the formula \eqref{equation: X61 t'} for
  $t'(\tau)$.
\end{proof}

In the next lemma we determine the ``modular equation'' of level $5$,
i.e., the polynomial relation between $t(\tau)$ and $t(\gamma\tau)$
for $\gamma\in\Gamma^\ast(\O)\SM1{-2}21\Gamma^\ast(\O)$.

\begin{Lemma} \label{lemma: X65->X61}
  The Shimura curve $X_6^\ast(5)$ has signature
  $(0;2^2,4^2)$. The ramification data of the covering $X_6^\ast(5)\to
  X_6^\ast(1)$ are as follows. \newline
\centerline{
\begin{tikzpicture}[scale=0.07]
  \node[place] (P6) [label=below:$P_6$, token=1] at (-70,0) {};
  \node[place] (P2) [label=below:$P_2$, token=1] at (-20,0) {};
  \node[place] (P4) [label=below:$P_4$, token=1] at (50,0) {};
  \node[place] (R6) [label=above:$R_6$, token=1] at (-70,20) {};
  \node[place] (Q21) [label=above:$Q_2$, token=1] at (-50,20) {};
  \node[place] (Q22) [label=above:$Q'_2$, token=1] at (-30,20)
  {};
  \node[place] (R21) [label=above:$R_2$, token=1] at (-10,20) {};
  \node[place] (R22) [label=above:$R'_2$, token=1] at (10,20) {};
  \node[place] (Q4) [label=above:$Q_4$, token=1] at (30,20) {};
  \node[place] (Q4') [label=above:$Q_4'$, token=1] at (50,20) {};
  \node[place] (R4) [label=above:$R_4$, token=1] at (70,20) {};
  \draw[line width=.6pt] (P6) --  node[left=2pt] {$6$} (R6);
  \draw[line width=.6pt] (Q21) --  node[left=2pt] {$1$}(P2);
  \draw[line width=.6pt] (Q22) --  node[left=2pt] {$1$}(P2);
  \draw[line width=.6pt] (R21) --  node[right=2pt] {$2$}(P2);
  \draw[line width=.6pt] (R22) --  node[right=2pt] {$2$}(P2);
  \draw[line width=.6pt] (Q4) -- node[left=2pt] {$1$} (P4);
  \draw[line width=.6pt] (Q4') -- node[left=2pt] {$1$} (P4);
  \draw[line width=.6pt] (R4) -- node[right=2pt] {$4$} (P4);
\end{tikzpicture}}
If we let $t$ the Hauptmodul on $X_6^\ast(1)$ with values
$0,-540,\infty$ at $P_6$, $P_2$, and $P_4$, respectively, and let
$u$ be a Hauptmodul on $X_6^\ast(5)$ with values $0$ and $\infty$ at
$R_6$ and $R_4$, respectively, then with a suitable scaling of $u$, we
have
$$
  t=\frac{(30u)^6}{1+18u+225u^2}.
$$
Moreover, the Atkin-Lehner involution $w_5$ switches the two elliptic
points $Q_2$ and $Q_2'$ of order $2$ and switches the two elliptic
points $Q_4$ and $Q_4'$ of order $4$, so that
$$
  w_5:u\longmapsto\frac{11u+2}{252u-11}.
$$
Finally, the polynomial relation between $t(\tau)$ and
$s(\tau)=t(w_5\tau)$ is given by the polynomial in Appendix A.
\end{Lemma}

\begin{proof} The Shimura curve $X_6(5)$ has totally
  $$
    \left(1-\JS{-4}2\right)\left(1-\JS{-4}3\right)
    \left(1+\JS{-4}5\right)=4
  $$
  CM points of discriminant $-4$. These are elliptic points of order
  $2$ on $X_6(5)$. The Atkin-Lehner involution $w_2$ fixes these
  points and the Atkin-Lehner involution $w_3$ switches them
  pairwise. Thus, $X_6^\ast(5)$ has $2$ elliptic points of order $4$.
  The curve $X_6(5)$ has no elliptic points of order $3$ since
  $\JS{-3}5=-1$. Thus, all the other elliptic points on $X_6^\ast(5)$
  are the fixed points of the Atkin-Lehner involutions $w_2$, $w_3$,
  and $w_6$, which, if exist, are CM points of discriminants $-8$, $-3$
  or $-12$, and $-24$, respectively. Since
  $\JS{-8}5=\JS{-3}5=\JS{-12}5=-1$, CM points of discriminant $-8$,
  $-3$, or $-12$ do not exist on $X_6(5)$. The number of CM points of
  discriminant $-24$ on $X_6(5)$ is
  $$
    2\left(1-\JS{-24}2\right)\left(1-\JS{-24}3\right)
     \left(1+\JS{-24}5\right)=4.
  $$
  (The integer $2$ stands for the class number of imaginary quadratic
  order of discriminant $-24$.)
  The Atkin-Lehner involution $w_6$ fixes these points, while the
  Atkin-Lehner involution $w_2$ switches them pairwise. Therefore,
  $X_6^\ast(5)$ has only $2$ elliptic points of order $2$ coming
  from CM points of discriminant $-24$. Then the genus formula shows
  that $X_6^\ast(5)$ has genus $0$. We conclude that $X_6^\ast(5)$ has
  signature $(0;2^2,4^2)$.

  The ramification data of $X_6^\ast(5)$ follow immediately from the
  above information.

  Now suppose that the Hauptmodul $t$ of $X_6^\ast(1)$ is chosen in a
  way that $t(P_6)=0$, $t(P_2)=-540$, and $t(P_4)=\infty$. If $u$ is a
  Hauptmodul on $X_6^\ast(5)$ with $u(R_6)=0$ and $u(R_4)=\infty$,
  then
  $$
     t=\frac{Au^6}{1+au+bu^2}
  $$
  for some complex numbers $A$, $a$, and $b$. Then the ramification
  data at $P_2$ imply that
  $$
    Au^6+540(1+au+bu^2)=540(1+cu+du^2)(1+eu+fu^2)^2
  $$
  for some complex numbers $c$, $d$, $e$, and $f$. To have nicer
  coefficients, we scale $u$ such that $a=18$. (The case $a=0$ yields
  $t=-270u^6/(1-3u^2/2)$, but then $w_5:u\to-u$, which implies that
  $t$ is a rational function of a Hauptmodul on $X_6^\ast(5)/w_5$. This
  is absurd.) Comparing the coefficients of the two sides above, we
  get $t=(30u)^6/(1+18u+225u^2)$.
\end{proof}

\begin{Lemma} \label{lemma: X61 eigenvalues}
  For $k=8,12,16,22,30,38$, let $f_k$ denote the
  automorphic form of weight $k$ on $\Gamma^\ast(\O)$ that spans the
  one-dimensional space $S_k(\Gamma^\ast(\O))$. Then the eigenvalues
  $\lambda$ of $T_5$ for $f_k$ are
  $$ \extrarowheight3pt
  \begin{array}{c|cccccc} \hline\hline
  k & 8 & 12 & 16 & 22 & 30 & 38 \\ \hline
  \lambda & -114 & 3630 & 77646 & -23245050
  & -21003872250 & 4477461318150 \\ \hline\hline
  \end{array}
  $$
\end{Lemma}

\begin{proof} By the Jacquet-Langlands correspondence (Proposition
  \ref{proposition: Jacquet-Langlands}),
  $$
    S_k(\Gamma^\ast(\O))\simeq S_k^\new(\Gamma_0(6),-1,-1),
  $$
  where $S_k^\new(\Gamma_0(6),-1,-1)$ denotes the Atkin-Lehner subspace of
  $S_k^\new(\Gamma_0(6))$ with eigenvalues $-1$ for both $W_2$ and
  $W_3$. We then look up the eigenvalues of $T_5$ in William Stein's
  modular form database \cite{Stein}. Alternatively, one can use the
  trace formulas of Eichler and Yamauchi
  \cite{Eichler-basis,Yamauchi-trace} to find the eigenvalues. (See
  \cite[Proposition 52]{Yang-correspondence} for a simplified trace
  formula, specifically for the group $\Gamma_0(6)$.)
\end{proof}

Using the informations above, we can determine which root of the
modular equation corresponds to $t(\gamma_j\tau)$.

\begin{Corollary} \label{corollary: X61 t}
  Let $\gamma_j=\SM{a_j}{b_j}{c_j}{d_j}$,
  $j=0,\ldots,5$, be the coset representatives given in Lemma
  \ref{lemma: X61 cosets}. There are rational numbers $A_n$,
  $n=0,1,2,\ldots$ with
  $$
    A_0=\frac{74649600}{14641}, \
    A_1=\frac{2918799360}{161051}, \
    A_2=\frac{69264688896}{1771561},\ \ldots,
  $$
  such that in a small neighborhood of $P_6$, the
  $t$-expansion of $t(\gamma_j\tau)$ is given by
  $$
    t(\gamma_j\tau)=\sum_{n=0}^\infty A_n(\zeta^jt^{1/6})^n,
  $$
  where $\zeta=e^{2\pi i/6}$ and $t^{1/6}$ is defined as in the
  paragraph preceding Lemma \ref{lemma: X61 t'}. In particular, at
  $\tau=\gamma_jP_6$, we have
  $t(\gamma_jP_6)=A_0=74649600/14641=2^{12}\cdot3^6\cdot5^2/11^4$.

  In addition, at $t=A_0$, we have
  \begin{equation} \label{equation: X61 F(t0)}
    F_1(A_0)-CF_2(A_0)=\sqrt[6]{\frac{11}{5^5}}(z+z^5-z^7), \qquad
    z=e^{2\pi i/24}.
  \end{equation}
\end{Corollary}

\begin{proof} Let $\Phi(s,t)$ be the modular equation given in
  Appendix A. We have
  $$
    \Phi(s,0)=-625(14641s-74649600)^6,
  $$
  which implies that $t(\gamma_jP_6)=74649600/14641=A_0$ for all $j$.
  Setting $s=\widetilde s+A_0$, the modular equation becomes
  \begin{equation} \label{equation: modular equation}
  \begin{split}
  &-625(11^6+O(t))^4\widetilde s^6+O(t)\widetilde
  s^5+O(t)\widetilde s^4+O(t)\widetilde s^3+O(t)\widetilde s^2
  +O(t)\widetilde s \\
  &\qquad
   +625t(2^{66}\cdot3^{36}\cdot5^6\cdot17^6\cdot23^6/11^6+O(t))=0.
  \end{split}
  \end{equation}
  Using Newton's polygon and Hensel's lemma, we see that each root of
  the equation $\widetilde s^6-A_1^6t=0$ lifts uniquely to a solution
  of \eqref{equation: modular equation}, where
  $A_1=2^{11}\cdot3^6\cdot5\cdot17\cdot23/11^5=2918799360/161051$. Since
  all coefficients in $\Phi(s,t)$ are rational numbers, the solution
  of \eqref{equation: modular equation} with the initial term
  $A_1t^{1/6}+\cdots$ have rational numbers as coefficients. We now
  show that the series $A_0+A_1t^{1/6}+A_2t^{2/6}+\cdots$ is the
  $t$-expansion of $t(\gamma_0\tau)$.

  By Theorem \ref{theorem: triangle}, the space
  $S_{12}(\Gamma^\ast(\O))$ is spanned by $F=(F_1-CF_2)^{12}$, where
  $F_1$, $F_2$, and $C$ are given as in Lemma \ref{lemma: X61 t'}.
  Now by Lemma \ref{lemma: X61 eigenvalues}, we have
  $$
    5^{11}\sum_{j=1}^6\frac1{(c_j\tau+d_j)^{12}}F(t(\gamma_j\tau))
   =3630F(t(\tau)),
  $$
  valid in a neighborhood of $P_6$, where $\gamma_j$ are the coset
  representatives given in Lemma \ref{lemma: X61 cosets}.
  Specializing $\tau$ to $P_6$ and using \eqref{equation: X61 ctau+d},
  we get
  $$
    5^{11}(z+z^5-z^7)^{-12}\sum_{j=1}^6F(t(\gamma_jP_6))=3630, \qquad
    z=e^{2\pi i/24}.
  $$
  Now we have $t(\gamma_jP_6)=A_0$ for all $j$. Thus,
  $$
    F(A_0)=\frac{11^2}{5^{10}}(z+z^5-z^7)^{12},
  $$
  i.e., $F_1(A_0)-CF_2(A_0)=z^{2j}(z+z^5-z^7)\sqrt[6]{11/5^5}$ for
  some $j$. Approximating numerically, we find this integer $j$ is
  equal to $0$. This proves \eqref{equation: X61 F(t0)}.

  Now assume that the $t$-expansion of $t(\gamma_j\tau)$ is
  $B_0+B_1t^{1/6}+\cdots$ with $B_0=A_0$. We have, by \eqref{equation: C},
  $$
    B_1=\lim_{\tau\to P_6}\frac{t(\gamma_0\tau)-B_0}
   {C^{-1}(\tau-P_6)/(\tau-\overline P_6)}.
  $$
  By L'Hopital rule and \eqref{equation: X61 ctau+d}, it is equal to
  \begin{equation} \label{equation: X61 identify temp}
    B_1=C(P_6-\overline P_6)\frac{5}{(c_0P_6+d_0)^2}t'(\gamma_0P_6).
  \end{equation}
  Combining \eqref{equation: X61 t'} and \eqref{equation: X61 F(t0)}, we
  find
  $$
    t'(\gamma_0P_6)=\frac{6(z+z^5-z^7)^2}{C(P_6-\overline P_6)}
    \sqrt[3]{\frac{11}{5^5}}B_0^{5/6}(1+B_0/540)^{1/2}
   =\frac{2^{11}\cdot3^6\cdot17\cdot23(z+z^5-z^7)^2}{11^5C(P_6-\overline
   P_6)}.
  $$
  Substituting this and \eqref{equation: X61 ctau+d} into
  \eqref{equation: X61 identify temp}, we arrive at
  $$
    B_1=\frac{2^{11}\cdot3^6\cdot5\cdot17\cdot23}{11^5}=A_1.
  $$
  This shows that the solution $A_0+A_1t^{1/6}+\cdots$ of the modular
  equation corresponds to $t(\gamma_0\tau)$. By \eqref{equation:
    gammaj} and \eqref{equation: sixth root}, it follows that
  $$
    t(\gamma_j\tau)=t(\gamma_0M_6^j\tau)
   =\sum_{n=0}^\infty A_nt(M_6^j\tau)^{n/6}
   =\sum_{n=0}^\infty A_n(\zeta^jt^{1/6})^{n}.
  $$
  This completes the proof.
\end{proof}

An interesting consequence of the above calculation is the following
evaluation of hypergeometric functions.

\begin{Corollary} \label{corollary: hypergeometric evaluation}
We have the evaluations
\begin{equation*}
\begin{split}
  _2F_1\left(\frac1{24},\frac7{24};\frac56;
  -\frac{2^{10}\cdot3^3\cdot5}{11^4}\right)&=\sqrt6
  \sqrt[6]{\frac{11}{5^5}}, \\
  {}_2F_1\left(\frac5{24},\frac{11}{24};\frac76;
  -\frac{2^{10}\cdot3^3\cdot5}{11^4}\right)&=
   \frac{\sqrt[3]2(1+\sqrt2)\sqrt[6]{11^5}}{20\sqrt3}
   \frac{\Gamma(7/6)\Gamma(13/24)\Gamma(19/24)}
   {\Gamma(5/6)\Gamma(17/24)\Gamma(23/24)}.
\end{split}
\end{equation*}
\end{Corollary}

\begin{proof} Let the notations $z$, $t_0$, $F_1$, and $F_2$ be the
  same as in the above corollary. By \eqref{equation: C}, we have
  $$
    \frac{CF_2(t_0)}{F_1(t_0)}=\frac{\gamma_jP_6-P_6}
   {\gamma_jP_6-\overline P_6}=\frac16(1+z^2+z^4-z^6).
  $$
  Combining this with \eqref{equation: X61 F(t0)} and simplifying, we
  get the claimed formulas.
\end{proof}

\begin{Remark} \label{remark: hypergeometric evaluation}
If we consider the Hecke operator $T_7$ instead, we will obtain analogous
formulas
\begin{equation*}
\begin{split}
  _2F_1\left(\frac1{24},\frac7{24};\frac56;
  \frac{2^{10}\cdot3^3\cdot5^6\cdot7}{11^4\cdot23^4}\right)&=
  \frac{2\sqrt2}7\sqrt[6]{11\cdot23}, \\
  {}_2F_1\left(\frac5{24},\frac{11}{24};\frac76;
  \frac{2^{10}\cdot3^3\cdot5^6\cdot7}{11^4\cdot23^4}\right)&=
   \frac{\sqrt[3]2(1+\sqrt2)}{140\sqrt3}
   \sqrt[6]{\frac{11^5\cdot23^5}7} \\
  &\qquad\qquad\times
   \frac{\Gamma(7/6)\Gamma(13/24)\Gamma(19/24)}
   {\Gamma(5/6)\Gamma(17/24)\Gamma(23/24)}.
\end{split}
\end{equation*}
We will not give a proof here.

Note that the numbers $-2^{10}\cdot3^3\cdot5/11^4$ and
$2^{10}\cdot3^3\cdot5^6\cdot7/(11^4\cdot23^4)$ correspond to the
CM-points of discriminants $-75$ and $-147$ on the Shimura curve
$X_6^\ast(1)$, respectively. In fact, with a little extra work, one
can show that at CM-points $\tau$ of discriminants $-3n^2$, $(n,6)=1$,
the values of ${}_2F_1(1/24,7/24;5/6;t(\tau)/540)$ are all algebraic
numbers. It will be an interesting problem to determine when $s$ and
${}_2F_1(1/24,7/24;5/6;s)$ are both algebraic over $\Q$.
\end{Remark}

The last information we need in order to compute Hecke operators is
the $t$-expansion of $F(\tau)=F_1(t(\gamma_j\tau))-CF_2(t(\gamma_j\tau))$ near
$P_6$. We will use the Jacquet-Langlands correspondence for this
purpose.

Set
$$
  f_0=F^{12}, \qquad f_1=t^{1/6}(1+t/540)^{1/2}F^{22}, \qquad
  f_2=t^{1/3}F^8,
$$
$$
  f_3=t^{1/2}(1+t/540)^{1/2}F^{30}, \qquad
  f_4=t^{2/3}F^{16}, \qquad f_5=t^{5/6}(1+t/540)^{1/2}F^{38}.
$$
By \eqref{equation: basis for X61 2}, these span the one-dimensional
spaces of automorphic forms of weights $12$, $22$, $8$, $30$, $16$,
and $38$, respectively. Let $k_\ell$ and $\lambda_\ell$,
$\ell=0,\ldots,5$, be the weights of $f_j$ and the eigenvalues for
$T_5$ given in Lemma \ref{lemma: X61 eigenvalues}. In other words, if
we let $\gamma_j=\SM{a_j}{b_j}{c_j}{d_j}$, $j=0,\ldots,5$ be the coset
representatives given in Lemma \ref{lemma: X61 cosets}, we have 
\begin{equation} \label{equation: recursion 1}
  5^{k_\ell/2-1}\sum_{j=0}^5\frac1{(c_j\tau+d_j)^{k_\ell}}f_\ell(\gamma_j\tau)
 =\lambda_\ell f_\ell(\tau)
\end{equation}
for $\ell=0,\ldots,5$. Note that the indices are arranged such that
\begin{equation} \label{equation: recursion 2}
  k_\ell/2+\ell\equiv 0\mod 6.
\end{equation}

Now the $t$-expansion of $\tau$ is known by
\eqref{equation: C}. Since $\gamma_jP_6$ lies on the boundary from
$P_6$ to $P_4$ of the fundamental domain, according to the agreement
on $t^{1/6}$ and $(1+t/540)^{1/2}$ made in the paragraph preceding
Lemma \ref{lemma: X61 t'}, $t(\gamma_jP_6)^{1/6}$ and
$(1+t(\gamma_jP_6)/540)^{1/2}$ are both positive. Then the
$t$-expansions of $t(\gamma_j\tau)^{1/6}$ and
$(1+t(\gamma_j\tau)/540)^{1/2}$ can be determined from that of
$t(\gamma_j\tau)$ given in Corollary \ref{corollary: X61 t}. They are
$$
  t(\gamma_j\tau)^{1/6}=12\sqrt[3]{\frac5{11^2}}\left(1
  +\frac{391}{660}\zeta^j t^{1/6}+\frac{14543}{36300}(\zeta^jt^{1/6})^2
  +\cdots\right)
$$
and
$$
  (1+t(\gamma_j\tau)/540)^{1/2}=\frac{391}{121}\left(1+
  \frac{6912}{4301}\zeta^jt^{1/6}+\frac{514656}{236555}(\zeta^jt^{1/6})^2
 +\cdots\right),
$$
respectively, where $\zeta=e^{2\pi i/6}$. Now assume that the
$t$-expansion of $F(\gamma_0\tau)$ near $P_6$ is
$B_0+B_1t^{1/6}+\cdots$. By \eqref{equation: gammaj} and
\eqref{equation: sixth root}, we have
$$
  F(\gamma_j\tau)=\sum_{n=0}^\infty B_n(\zeta^jt^{1/6})^n.
$$
We now determine $B_n$ inductively.

The value of $B_0$ is already determined in Corollary \ref{corollary:
X61 t}. It is equal to $(z+z^5-z^7)\sqrt[6]{11/5^5}$, where $z=e^{2\pi
  i/24}$. Now assume that the values of $B_m$ are known up to
$m=n-1$. To determine $B_n$, we let $\ell\in\{0,\ldots,5\}$ be the
integer satisfying $\ell\equiv n\mod 6$ and consider \eqref{equation:
  recursion 1}. The coefficient of $t^{n/6}$ on the left-hand side of
\eqref{equation: recursion 1} is equal to
\begin{equation} \label{equation: recursion 3}
  (\text{a known number})+
  \sum_{j=0}^5\frac{5^{k_\ell/2-1}k_\ell}{(c_jP_6+d_j)^{k_\ell}}
  C_\ell B_0^{k_\ell-1}B_n\zeta^{jn},
\end{equation}
where
$$
  C_\ell=12^\ell(5/121)^{\ell/3}\times  \begin{cases}
  1, &\text{if }\ell\equiv 0\mod 2, \\
  391/121, &\text{if }\ell\equiv 1\mod 2. \end{cases}
$$
By \eqref{equation: X61 ctau+d}, we have
$$
  (c_jP_6+d_j)^{-k_\ell}=\zeta^{jk_\ell/2}(c_0P_6+d_0)^{-k_\ell}.
$$
In view of \eqref{equation: recursion 2}, \eqref{equation: recursion
  3} is equal to
$$
  (\text{a known number})+\frac{6\cdot5^{k_\ell/2-1}k_\ell}{(c_0P_6+d_0)^{k_\ell}}
  C_\ell B_0^{k_\ell-1}B_n.
$$
This number must be equal to the coefficient of $t^{n/6}$ on the
right-hand side of \eqref{equation: recursion 1}. This determines the
value of $B_n$ inductively. The first few $B_n$ are given in Appendix B.

In general, if we wish to compute the Hecke operator $T_5$ on the
space of automorphic forms of weight $k$ on $X_6^\ast(1)$ with
dimension $d_k$, we just have to determine the $t$-expansions of $\tau$,
$t(\gamma_j\tau)$ and $F(\gamma_j\tau)$ up to the term
$t^{d_k-1+\{5k/12\}}$ and then express
$$
  5^{k/2-1}\sum_{j=0}^5\frac1{(c_j\tau+d)^k}g_\ell
$$
as a linear combination of $g_m$ by comparing the coefficients up to
the term $t^{d_k-1+\{5k/12\}}$ for each $g_\ell$ in \eqref{equation:
  basis for X61 2}. In Appendix C, we give the matrices for $T_5$ up
to weight $48$.

\begin{Remark} Using the Jacquet-Langlands correspondence, it is easy
  to deduce the matrices for other Hecke operators from that of $T_5$.
  For example, for the case of weight $24$, according to William
  Stein's modular form database \cite{Stein}, the pair of
  Galois-conjugate normalized Hecke eigenforms in
  $S_{24}^\new(6,-1,-1)$ have Fourier expansions
  $$
    q + 2048q^2 + 177147q^3 + \cdots + aq^5 + \cdots + (-25a +
    3197833334)q^7 + \cdots,
  $$
  where $a$ is a root of the characteristic polynomial of $T_5$, which
  is irreducible over $\Q$. Now, the matrix for $T_5$ relative to our
  basis of automorphic forms on $X_6^\ast(1)$ is
  $$
    A=\begin{pmatrix}
   10980750 & 3111696/5 \\ 55987200000 &  14267406\end{pmatrix}.
  $$
  Thus, the matrix for $T_7$ relative to the same basis is
  $$
    -25A+3197833334=\begin{pmatrix}
    2923314584 & -15558480 \\ -1399680000000 & 2841148184 \end{pmatrix}.
  $$
\end{Remark}
\end{section}

\begin{section}{Hecke operators on $X_{10}^\ast(1)$}
  In this section, we will consider the Shimura curve
  $X^\ast_{10}(1)$. The argument runs completely parallel to the case
  of $X^\ast_6(1)$, so we will just sketch our computation.

  As in \cite[Section 5.5.3]{Alsina-Bayer}, we let $B$ be the algebra
  generated by $I$ and $J$ over $\Q$ with the relations
  $$
    I^2=2, \qquad J^2=5, \qquad IJ=-JI.
  $$
  Then $B$ is a quaternion algebra of discriminant $10$ over $\Q$.
  Fix the maximal order $\O$ to be $\Z+\Z I+\Z(1+J)/2+\Z(I+IJ)/2$ and
  choose the embedding $\iota:B\to M(2,\R)$ to be
  $$
    I\longmapsto\M{\sqrt2}00{-\sqrt2}, \qquad
    J\longmapsto\M0150.
  $$
  A fundamental domain for $\Gamma(\O)$ is given in
  \cite[Section 5.5.3]{Alsina-Bayer}, from which we deduce that a
  fundamental domain for $\Gamma^\ast(\O)$ is

  \centerline{\epsfig{file=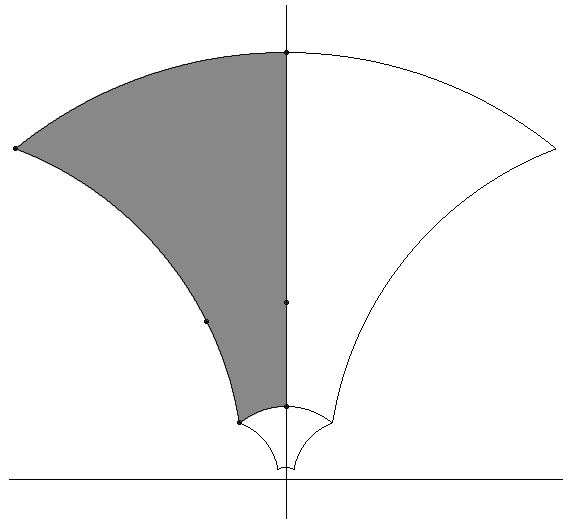,height=2.7in,width=3in}}

  \noindent Here the grey area represents a fundamental domain for
  $X_{10}^\ast(1)$. The six marked points on the boundary, in
  clockwise order from the top one on the imaginary axis, are
  $$
    \frac{(\sqrt2+1)i}{\sqrt5}, \quad \frac i{\sqrt 5}, \quad
    \frac{(\sqrt2-1)i}{\sqrt5}, \quad
    \frac{-\sqrt2+\sqrt 3i}{5(\sqrt2+1)}, \quad
    \frac{-1+2i}5, \quad\frac{-\sqrt2+\sqrt3i}{5(\sqrt2-1)},
  $$
  respectively. The grey area and the three other white areas form a
  fundamental domain for $X_{10}(1)$.

  The representatives of the elliptic point of order $3$ and the
  three elliptic points of order $2$ are
  $$
    P_3=\frac{-\sqrt2+\sqrt 3i}{5(\sqrt2-1)}, \qquad
    P_2=\frac{-1+2i}5, \qquad P_2'=\frac i{\sqrt5(\sqrt2-1)}, \qquad
    P_2''=\frac i{\sqrt 5}
  $$
  with the isotropy subgroups generated by
  \begin{equation} \label{equation: 10 isotropy}
  \begin{split}
    M_3=\frac12\M{-1-\sqrt2}{-1-\sqrt2}{5(-1+\sqrt2)}{-1+\sqrt2}, \qquad
   &M_2=\frac1{2\sqrt2}\M{\sqrt2}{\sqrt2}{-5\sqrt2}{-\sqrt2}, \\
   M_2'=\frac1{\sqrt5}\M0{1-\sqrt2}{5(1+\sqrt2)}0, \qquad
  &M_2''=\frac1{\sqrt{10}}\M0{\sqrt2}{-5\sqrt2}0,
  \end{split}
  \end{equation}
  respectively. Note that the points $P_2$, $P_2'$, and $P_2''$ are the
  fixed points of the Atkin-Lehner involutions $w_2$, $w_5$, and
  $w_{10}$, respectively. That is, they are CM-points of discriminants
  $-8$, $-20$, and $-40$, respectively. According to \cite{Elkies},
  there is a Hauptmodul $t(\tau)$ on $\Gamma^\ast(\O)$ that takes
  values $0$, $\infty$, $2$, and $27$ at $P_3$, $P_2$, $P_2'$, and
  $P''$, respectively. Also, by \eqref{equation: Q10}, the Schwarzian
  differential equation associated to $t$ is
  $$
    \frac{d^2}{dt^2}f+\frac{3t^4-119t^3+3157t^2-7296t+10368}
   {16t^2(t-2)^2(t-27)^2}f=0.
  $$
  In other words, near the point $P_3$, the $t$-expansion of
  $t'(\tau)$ is the square of a linear combination of two solutions
  \begin{equation*}
  \begin{split}
    F_1(t)&=t^{1/3}\left(1-\frac{10}{81}t-\frac{18539}{839808}t^2
   -\frac{168605}{25509168}t^3-\frac{107269219465}{46548313473024}t^4
   +\cdots\right), \\
    F_2(t)&=t^{2/3}\left(1-\frac5{81}t-\frac{99095}{5878656}t^2
   -\frac{8353325}{1428513408}t^3-\frac{851170821485}{385081502367744}t^4
   +\cdots\right)
  \end{split}
  \end{equation*}
  of the differential equation above. To determine the linear
  combination, we follow the computation in Lemma \ref{lemma: X61 t'}.

  Similar to \eqref{equation: C}, we have
  $$
    \frac{\tau-P_3}{\tau-\overline P_3}=C\frac{F_2}{F_1},
  $$
  where $C$ is a nonzero constant. Thus,
  $$
    \frac{d\tau}{dt}=C(P_3-\overline P_3)
    \frac{F_1dF_2/dt-F_2dF_1/dt}{(F_1-CF_2)^2}.
  $$
  Here because the differential equation is normalized, the numerator
  $F_1dF_2/dt-F_2dF_1/dt$ is just a constant. In fact, by computing
  the leading coefficients, we find that it is $1/3$. Thus,
  \begin{equation} \label{equation: 10 t'}
    t'(\tau)=\frac{3(F_1-CF_2)^2}{C(P_3-\overline P_3)}.
  \end{equation}
  Note that the function $t(\tau)^{1/3}$ is defined in a way such that
  it takes negative values along the arc from $P_3$ to $P_2$ and
  becomes a holomorphic function near the point $P_3$.

  Now, by Theorem \ref{theorem: basis},
  for an even integer $k\ge 4$, a basis for $S_k(\Gamma^\ast(\O))$ is
  \begin{equation} \label{equation: 10 basis}
    \frac{t^j\left(F_1(t)-CF_2(t)\right)^{k/2}}
   {t^{\FL{k/3}}(1-t/2)^{\FL{k/4}}(1-t/27)^{\FL{k/4}}}, \qquad j=0,\ldots, d_k-1,
  \end{equation}
  for some constant $C$, where $d_k=1-k+\FL{k/3}+3\FL{k/4}$ is the
  dimension of $S_k(\Gamma^\ast(\O))$. We now compute the Hecke
  operator $T_3$ with respect to this basis.

  Let $\gamma=1+2I+IJ$. An Eichler order of level $3$ is given by
  $\O\bigcap\gamma^{-1}\O\gamma$. Choose the coset representatives of
  $\Gamma^\ast(\O)\backslash\Gamma^\ast(\O)\SM{1+2\sqrt2}{\sqrt2}
  {-5\sqrt2}{1-2\sqrt2}\Gamma^\ast(\O)$ to be
  \begin{equation*}
  \begin{split}
  \gamma_0=\frac12\M{3+\sqrt2}{1+\sqrt2}{5-5\sqrt2}{3-\sqrt2},
  \quad &
  \gamma_1=\frac1{\sqrt{10}}\M{-5}{-1-\sqrt2}{-5+5\sqrt2}{-5}, \\
  \gamma_2=\frac1{2\sqrt{10}}
  \M{5\sqrt2}{4+5\sqrt2}{20-25\sqrt2}{-5\sqrt2}, \quad &
  \gamma_3=\frac1{2\sqrt{10}}
  \M{10-5\sqrt2}{-2-3\sqrt2}{-10+15\sqrt2}{10+5\sqrt2}.
  \end{split}
  \end{equation*}
  These coset representatives have the properties
  $$
    \gamma_0P_3=P_3, \qquad
    \gamma_1P_3=\gamma_2P_3=\gamma_3P_3=\frac{-2\sqrt2+3\sqrt3i}
    {5(2\sqrt2-1)},
  $$
  and
  $$
    \gamma_2=\gamma_1M_3, \qquad \gamma_3=\gamma_1M_3^2,
  $$
  where $M_3$ is the generator of the isotropy subgroup of $P_3$ given
  in \eqref{equation: 10 isotropy}. Also, for $j=0,\ldots,3$, if we
  write $\gamma_j=\SM{a_j}{b_j}{c_j}{d_j}$, then
  \begin{equation} \label{equation: 10 cP3+d}
    c_0P_3+d_0=\frac{3-\sqrt3i}2, \qquad
    c_jP_3+d_j=\zeta^{j-1}\frac{-5-\sqrt2+\sqrt3i}{\sqrt{10}}, \qquad
    \zeta=e^{2\pi i/3}.
  \end{equation}

  By \eqref{equation: 10 t-u} and \eqref{equation: 10 w3}, the modular
  equation of level $3$ is the relation between
  $$
    t=\frac{216(u-1)^3}{(u+1)^2(9u^2-10u+17)}, \qquad
    s=\frac{216(1/9-u)^3}{(19/9-u)^2(9u^2-10u+17)}.
  $$
  We find that it is equal to
  \begin{equation*}
  \begin{split}
  &t(25t+192)^3+(7077888+2908160t-12612480t^2+1674720t^3-36750t^4)s
   \\
  &\qquad+(2764800-12612480t+8025390t^2-798210t^3+21609t^4)s^2 \\
  &\qquad+(360000+1674720t-798210t^2+33614t^3)s^3+(147t-125)^2s^4=0.
  \end{split}
  \end{equation*}
  As an equation in $s$, the $4$ roots of the above equation are
  \begin{equation} \label{equation: s0}
  s_0=-t - \frac{10}{27}t^2 - \frac{100}{729}t^3 -
  \frac{16675}{314928}t^4 - \frac{90125}{4251528}t^5-\cdots
  \end{equation}
  and
  \begin{equation} \label{equation: sj}
    s_j=-\frac{192}{25}-\frac{2992}{125}\zeta^{j-1}t^{1/3}
    -\frac{25044}{625}(\zeta^{j-1}t^{1/3})^2
    -\frac{501163}{9375}(\zeta^{j-1}t^{1/3})^3-\cdots
  \end{equation}
  for $j=1,2,3$. It is clear that $s_0$ is the $t$-expansion of
  $t(\gamma_0\tau)$ near $P_3$. To determine how $s_j$ are matched
  with $t(\gamma_m\tau)$ for $j,m=1,2,3$, we consider the space of
  automorphic forms on $X_{10}^\ast(1)$ of weight $18$. The space
  $S_{18}(\Gamma(\O))$ has dimension $1$ and is spanned by
  $$
    f_{18}=\frac{(F_1(t)-CF_2(t))^{18}}{t^6(1-t/2)^4(1-t/27)^4}
  $$
  By the Jacquet-Langlands correspondence and the modular form
  database of William Stein's \cite{Stein}, we have
  $T_3f_{18}=-14976f_{18}$. In other words,
  $$
    3^{17}\sum_{j=0}^3\frac1{(c_j\tau+d_j)^{18}}f_{18}(\gamma_j\tau)
   =-14976f_{18}(\tau).
  $$
  Evaluating the two sides at $\tau=P_3$ and using \eqref{equation: 10
    cP3+d}, we get
  $$
    3^{17}\left(-\frac1{3^9}
   +\frac{3\cdot10^9(F_1(-192/25)-CF_2(-192/25))^{18}}
    {(5+\sqrt2-\sqrt3i)^{18}(192/25)^6(121/25)^4(289/225)^4}\right)
   =-14976.
  $$
  Therefore,
  \begin{equation} \label{equation: 10 t'(P3)}
    (F_1(-192/25)-CF_2(-192/25))^2=-\epsilon
    \frac{2^3\cdot11\cdot17\cdot(5+\sqrt2-\sqrt3i)^2}{3^2\cdot5^4}
  \end{equation}
  for some $9$th root of unity $\epsilon$. Approximating numerically,
  we find this root of unity is equal to $1$. If we assume that the
  $t$-expansion of $t(\gamma_j\tau)$, $j=1,2,3$, is
  $B_0+B_1t^{1/3}+\cdots$, then, similar to \eqref{equation: X61
    identify temp},
  $$
    B_1=C(P_3-\overline P_3)\frac3{(c_jP_3+d_j)^2}t'(\gamma_jP_3).
  $$
  From \eqref{equation: 10 cP3+d}, \eqref{equation: 10 t'} and
  \eqref{equation: 10 t'(P3)}, it follows that
  $$
    B_1=-\zeta^{j-1}\frac{2992}{125}, \qquad \zeta=e^{2\pi i/3}.
  $$
  Therefore, we have $t(\gamma_j\tau)=s_j$, where $s_j$ is the power
  series in \eqref{equation: sj}.

  The last piece of information needed for computing the Hecke
  operator $T_3$ is the $t$-expansion of $F(t(\gamma_j\tau))$ near
  $P_3$, where $F(t)=F_1(t)-CF_2(t)$. For $\gamma_0$, since
  $t(\gamma_0P_3)=0$, the expansion is just
  $F_1(s_0)-CF_2(s_0)$, where $s_0$ is given by \eqref{equation: s0}.
  For $\gamma_1,\gamma_2,\gamma_3$, we use the Jacquet-Langlands
  correspondence. Set
  $$
    f_0=\frac{F^{18}}{t^6(1-t/2)^4(1-t/27)^4}, \quad
    f_1=\frac{F^4}{t(1-t/2)(1-t/27)}, \quad
    f_2=f_1^2.
  $$
  These functions span the one-dimensional spaces of automorphic forms
  of weights $18$, $4$, and $8$, respectively. By the
  Jacquet-Langlands correspondence and the modular form database, we
  have
  $$
    T_3f_0=-14976f_0, \qquad T_3f_1=-8f_1, \qquad
    T_3f_2=28f_2.
  $$
  Using the same idea in the case of $X_6^\ast(1)$, we can inductively
  determine the $t$-expansion of $F(t(\gamma_j\tau))$. This is enough
  to compute $T_3$ for automorphic forms of general weights. The
  matrices for weights up to $32$ are given in Appendix D.

\begin{Remark} Similar to the case of $X_6^\ast(1)$, the Hecke
  operator $T_3$ gives rise to special values of $F_1(t)$ and $F_2(t)$
  at $t=-192/25$. However, because there does not seem to be a simple
  description of $F_1$ and $F_2$, we do not work out the values here.
\end{Remark}
\end{section}

\bigskip
\centerline{\sc Appendix A. Modular equation of level 5 for $X_6^\ast(1)$}
\bigskip

The relation between $t(\tau)$ and $s(\tau)=t(w_5\tau)$ in Lemma
\ref{lemma: X65->X61} is given by
$\Phi(s,t)=a_0(t)+a_1(t)s+\cdots+a_6(t)s^6=0$ with
$$
  a_0(t)=625(14641t-74649600)^6,
$$
\begin{equation*}
\begin{split}
  a_1(t)=&-127273923718594838908526411749785600000000000000 \\
 &-80334423751172973765172796218289002905600000000t\\
 &-429964557500791635545687183398954598400000000t^2\\
 &-584929357876511069442449306458521600000000t^3\\
 &-157301324859052802277036509414400000000t^4\\
 &-6539287187545403426159665668843750t^5\\
 &-10812982790452826706563610000t^6,
\end{split}
\end{equation*}
\begin{equation*}
\begin{split}
  a_2(t)=&62405475620899075027184847421440000000000000\\
  &-429964557500791635545687183398954598400000000t\\
  &+1785753116574594713648207726896742400000000t^2\\
  &+14568151384869872301210980765577600000000t^3\\
  &+16627948765574028094821145899437109375t^4\\
  &-4345461774128783231852293307250000t^5\\
  &+7122260106437394560116860000t^6,
\end{split}
\end{equation*}
\begin{equation*}
\begin{split}
  a_3(t)=&-16319418877271650617873609523200000000000\\
  &-584929357876511069442449306458521600000000t\\
  &+14568151384869872301210980765577600000000t^2\\
  &+6765662887332547989803108862517187500t^3\\
  &-126311317151610531211635483472500000t^4\\
  &-151264183450476527706794072400000t^5\\
  &-2085007542703940658656160000t^6,
\end{split}
\end{equation*}
\begin{equation*}
\begin{split}
  a_4(t)=&2400541447463893678227554304000000000\\
  &-157301324859052802277036509414400000000t\\
  &+16627948765574028094821145899437109375t^2\\
  &-126311317151610531211635483472500000t^3\\
  &-224960987576019075545123651160000t^4\\
  &+168218287650857617198656825600t^5\\
  &+228890990438717652531360000t^6,
\end{split}
\end{equation*}
\begin{equation*}
\begin{split}
  a_5(t)=&-188326942581441118735691136000000\\
  &-6539287187545403426159665668843750t\\
  &-4345461774128783231852293307250000t^2\\
  &-151264183450476527706794072400000t^3\\
  &+168218287650857617198656825600t^4\\
  &+386001766275853449228885504t^5,
\end{split}
\end{equation*}
and
$$
  a_6(t)=625(777924t-1771561)^4.
$$
\bigskip

\centerline{\sc Appendix B. $t$-expansion of automorphic forms on
  $X_6^\ast(1)$}
\bigskip

Let $z=e^{2\pi i/24}$. The first few terms in
\begin{equation*}
\begin{split}
  F(\tau):=F_1(t(\gamma_0\tau))-CF_2(\gamma_0\tau)=
  \sqrt[6]{\frac{11}{5^5}}(z+z^5-z^7)\left(1+
  \sum_{n=1}^\infty a_nt^{n/6}\right)
\end{split}
\end{equation*}
are
\begin{equation*}
\begin{split}
  a_1&=\frac C{13}(8z^6 - 3z^4 - 2z^2 - 9) - \frac7{55} \\
  a_2&=\frac C{715}(-56z^6 + 21z^4 + 14z^2 + 63) - \frac{161}{3025} \\
  a_3&=\frac C{39325}(-1288z^6 + 483z^4 + 322z^2 + 1449) - \frac{1379}{45375} \\
  a_4&=\frac C{589875}(-11032z^6 + 4137z^4 + 2758z^2 + 12411) -
       \frac{118027}{6655000} \\
  a_5&=\frac C{6655000}(-72632z^6 + 27237z^4 + 18158z^2 + 81711)
      -\frac{25165}{2108304} \\
  a_6&=\frac C{27407952}(-201320z^6 + 75495z^4 + 50330z^2 + 226485) - 
    \frac{219273755477}{26090262000000}.
\end{split}
\end{equation*}
\bigskip

\centerline{\sc Appendix C. Matrices for $T_5$ on $X_6^\ast(1)$ up to
  weight $48$}
\bigskip

  Here we list the matrices for $T_5$ on $X_6^\ast(1)$ up to weight
  $48$, computed using the recipe described in Section \ref{section:
    246}.

  Let $d_k=1-k+\lfloor k/4\rfloor+\lfloor 3k/8\rfloor
  +\lfloor 5k/12\rfloor$ be the dimension of the space of automorphic
  forms of weight $k$ on $X_6^\ast(1)$ and $g_\ell$,
  $\ell=0,\ldots,d_k-1$ be the basis given in \eqref{equation: basis
    for X61 2}. The matrices listed in the table below satisfy
  $$
    T_5\begin{pmatrix}g_0\\ \vdots\\ g_{d_k-1}\end{pmatrix}
  =M\begin{pmatrix}g_0\\ \vdots\\ g_{d_k-1}\end{pmatrix}.
  $$
$$ \extrarowheight3pt
\begin{array}{c|l} \hline\hline
k & M \\ \hline
8 & -114 \\
12& 3630 \\
16& 77646 \\
20& 1953390 \\
22& - 23245050\\
24& \displaystyle\begin{pmatrix}
       10980750 &  3111696/5 \\
    55987200000 &   14267406 \end{pmatrix} \\
28& 1220703150 \\
30& - 21003872250 \\
32& \displaystyle\begin{pmatrix}
    105068988750 &  376515216/5 \\
  12317184000000 & -39127734834 \end{pmatrix} \\
34& - 249151856250 \\
36& \displaystyle\begin{pmatrix}
      33216768750 & 44743076784/5 \\
  169361280000000 & 2408347964910 \end{pmatrix} \\
38& 4477461318150 \\
40& \displaystyle\begin{pmatrix}
    -70619784011250  &  45558341136/5 \\
  30930128640000000  &  36422537206926 \end{pmatrix}\\
42& 9372398943750 \\
44& \displaystyle\begin{pmatrix}
     896721261768750 &  -10018108383696/5 \\
 -321103388160000000 &  -695085225669330 \end{pmatrix} \\
46& \displaystyle\begin{pmatrix}
   -1294661994656250  &    -3322118218608 \\
 -181517500800000000  & -5089104194777850 \end{pmatrix} \\
48& \displaystyle\begin{pmatrix}
       100480725468750 &       225950546273760  &      2420662999104/5 \\
    512317872000000000 &     22159766272716750  &      5512559277456/5 \\
2612138803200000000000 &  -7950573190656000000  &   -23013714467131314
  \end{pmatrix} \\
\hline\hline
\end{array}
$$
\bigskip

\centerline{\sc Appendix D. Matrices for $T_3$ on $X_{10}^\ast(1)$ up to
  weight $32$}
\bigskip

  Here we list the matrices for $T_3$ on $X_{10}^\ast(1)$ up to weight
  $32$. Let $d_k=1-k+\lfloor k/3\rfloor+3\lfloor k/4\rfloor$ be the
  dimension of the space of automorphic forms of weight $k$ on
  $X_{10}^\ast(1)$ and $g_\ell$, $\ell=0,\ldots,d_k-1$ be the basis
  given in \eqref{equation: 10 basis}. The matrices listed in the
  table below satisfy
  $$
    T_3\begin{pmatrix}g_0\\ \vdots\\ g_{d_k-1}\end{pmatrix}
  =M\begin{pmatrix}g_0\\ \vdots\\ g_{d_k-1}\end{pmatrix}.
  $$
$$ \extrarowheight3pt
\begin{array}{c|l} \hline\hline
k & M \\ \hline
4 & -8 \\
8 & 28 \\
12& \begin{pmatrix}
  468  & -98 \\
 -1728 & 136 \end{pmatrix} \\
16& \begin{pmatrix}
  1728 &  490 \\
 34560 &-3572 \end{pmatrix} \\
18& -14976 \\
20& \begin{pmatrix}
   -2268 &  -2450 \\
 -328320 &  35992 \end{pmatrix} \\
22& -21924 \\
24& \begin{pmatrix}
  227772 & -272244 &   14406 \\
 -388800 & -258192 &   12250 \\
 2985984 &  711936 & -199556 \end{pmatrix}\\
26& 162864 \\
28& \begin{pmatrix}
     420552 &   949620 &  -72030 \\
    -933120 &  4479732 &  -61250 \\
 -104509440 & 31147200 & -196568 \end{pmatrix} \\
30& \begin{pmatrix}
 -6676344 & 4593750 \\
 14541120 & 3031596 \end{pmatrix} \\
32& \begin{pmatrix}
   29821932 &   -5456052 &   360150 \\
   95084928 &  -48253536 &   306250 \\
 1803534336 & -618444288 & 19290988 \end{pmatrix} \\ \hline\hline
\end{array}
$$

\end{document}